\documentclass[onefignum,onetabnum]{siamart190516}

\usepackage{lipsum}
\usepackage{amsfonts}
\usepackage{graphicx}
\usepackage{epstopdf}
\ifpdf
  \DeclareGraphicsExtensions{.eps,.pdf,.png,.jpg}
\else
  \DeclareGraphicsExtensions{.eps}
\fi

\usepackage{amssymb}
\usepackage{amsfonts}
\usepackage{bm}
\usepackage{bbm}
\usepackage{xcolor}
\usepackage{enumitem}

\usepackage{xspace}

\usepackage[font={small,it}]{caption}
\usepackage{algorithm}
\usepackage[noend]{algpseudocode}

\newsiamremark{remark}{Remark}
\newsiamthm{claim}{Claim}

\usepackage{tabularx}
\newcolumntype{C}{>{\centering\arraybackslash}p{19mm}}
\newcolumntype{G}{>{\centering\arraybackslash}p{4mm}}
\newcolumntype{S}{>{\centering\arraybackslash\scriptsize}p{4mm}}

\usepackage{multirow}

\headers{Accessibility and Controllability of Statistical Linearization}{R. Bonalli, C. Leparoux, B. H\'eriss\'e, and F. Jean}

\title{On the Accessibility and Controllability of Statistical Linearization for Stochastic Control: Algebraic Rank Conditions and their Genericity}

\author{R. Bonalli\thanks{Laboratoire des Signaux et Syst\`emes, Université Paris-Saclay, CNRS, CentraleSupélec, Bât. Bréguet, 3 Rue Joliot Curie, 91190 Gif-sur-Yvette, France. \textit{E-mail}: {\tt riccardo.bonalli@cnrs.fr}.} \and C. Leparoux\thanks{UMA, ENSTA Paris, Institut Polytechnique de Paris, 91120 Palaiseau, France; and DTIS, ONERA, Université Paris-Saclay, 91123 Palaiseau, France. \textit{E-mail}: {\tt clara.leparoux@ensta-paris.fr}.} \and B. H\'eriss\'e\thanks{DTIS, ONERA, Université Paris-Saclay, F-91123 Palaiseau, France. \textit{E-mail}: {\tt bruno.herisse@onera.fr}.} \and F. Jean\thanks{UMA, ENSTA Paris,
        Institut Polytechnique de Paris, 91120 Palaiseau, France. \textit{E-mail}: {\tt frederic.jean@ensta-paris.fr}.}}

\usepackage{amsopn}

\usepackage{lipsum}
\newcommand\blfootnote[1]{%
  \begingroup
  \renewcommand\thefootnote{}\footnote{#1}%
  \addtocounter{footnote}{-1}%
  \endgroup
}
\RequirePackage{lineno}

\usepackage{enumitem}
\newenvironment{ItemizE}{\begin{itemize}[itemsep=2.5pt,topsep=5pt,left=0ex]}{\end{itemize}}

\usepackage{todonotes}  

\reversemarginpar

\usepackage{appendix}

\newcommand{\R}{\mathbb{R}}
\newcommand{\Sym}{\textnormal{Sym}}
\newcommand{\Skew}{\textnormal{Skew}}
\newcommand{\U}{\mathcal{U}}
\newcommand\A{\mathcal{A}}

\newcommand\lie{\mathrm{Lie}}
\newcommand\I{\mathcal{I}}
\newcommand\F{\mathfrak{F}}
\newcommand\G{g}
\newcommand\M{V}

\newcommand\pder[2]{\frac{\partial #1}{\partial #2}}
\newcommand\vecdef[2]{\begin{pmatrix} #1 \\ #2 \end{pmatrix}}
\newcommand\Vect[1]{\mathrm{span}\left\lbrace #1\right\rbrace}
\DeclareMathAlphabet{\mathpzc}{OT1}{pzc}{m}{it}
\newcommand\f{\mathfrak{f}}

\newcommand{\codim}{\textnormal{codim} \,}
\newcommand{\ad}{\textnormal{ad}}
\newcommand{\rank}{\textnormal{rank} \,}
\newcommand{\ncontrol}{n_u} 
\newcommand{\nWiener}{d}
\newcommand{\Id}{I}

\theoremstyle{plain}

\ifpdf
\hypersetup{
  pdftitle=On the Accessibility and Controllability of Statistical Linearization for Stochastic Control: Algebraic Rank Conditions and their Genericity,
  pdfauthor={R. Bonalli, C. Leparoux, B. H\'eriss\'e, and F. Jean}
}
\fi

\begin{document}

\maketitle
\blfootnote{This paper was submitted for review on \today.}

\begin{abstract}
Statistical linearization has recently seen a particular surge of interest as a numerically cheap method for robust control of stochastic differential equations. Although it has already been successfully applied to control complex stochastic systems, accessibility and controllability properties of statistical linearization, which are key to make the robust control problem well-posed, have not been investigated yet. In this paper, we bridge this gap by providing sufficient conditions for the accessibility and controllability of statistical linearization. Specifically, we establish simple sufficient algebraic conditions for the accessibility and controllability of statistical linearization, which involve the rank of the Lie algebra generated by the drift only. In addition, we show these latter algebraic conditions are essentially sharp, by means of a counterexample, and that they are generic with respect to the drift and the initial condition.
\end{abstract}

\begin{keywords}
    Robust stochastic control, statistical linearization, accessibility and controllability, algebraic rank conditions, generic results.
\end{keywords}

\begin{AMS}
  93B03, 93B27, 93E03.
\end{AMS}


\section{Introduction}
Let $(\Omega,\mathcal{F},(\mathcal{F}_t)_{t \ge 0},\mathbb{P})$ be a filtered probability space, with complete filtration $(\mathcal{F}_t)_{t \ge 0}$, and let $(W_t)_{t \ge 0}$ be a $\nWiener$-dimensional Wiener process which is adapted to $(\mathcal{F}_t)_{t \ge 0}$. In this work, we focus on stochastic control systems which are modeled by It\^o-type stochastic differential equations of the form
\begin{equation} \label{eq:eds}
\mathrm{d}x(t) = f(x(t),u(t)) \; \mathrm{d}t + \G(x(t)) \; \mathrm{d}W_t .
\end{equation}
where 
$f$ and $\G$ denote drift and diffusion, respectively, whereas $x$ and $u$ denote the finite-dimensional state of the system and the control applied to it, respectively.

Equations such as \eqref{eq:eds} are frequently leveraged to accurately model control systems which are subject to exogenous random perturbations, finding beneficial use in many applications which include aerospace \cite{Ridderhof2019}, biology \cite{Berret20201}, and finance \cite{Carmona2016}, to name a few. In particular, one may often mitigate harmful random outcomes which perturb control systems by appropriately controlling \eqref{eq:eds} with the objective of minimizing (or, at least constraining) the variance of the state $x$ throughout the control time interval: we call this latter \textit{minimal variance control}. Successful application of this technique has been originally obtained in the context of linear systems, enabling numerical approaches ranging from exact covariance control \cite{Hotz1985} and variance minimization through linear-quadratic regulators \cite{Sain1966,Chen2015}, to multi-objective variance minimization via sensitive robustness \cite{Nagy2004}. For what concerns non-linear systems, we find extension of this latter paradigm \cite{Rustem1994}, followed by open-loop-based sensitive robustness \cite{Darlington2000}, and more recently, by improved perturbation-based methods which penalize an appropriate sensitivity matrix within the cost \cite{Seywald2019}.

Unfortunately, the aforementioned procedures suffer from expensive numerical computations, which has recently fostered the investigation of methods for minimal variance control which rely on numerically cheaper approximations of the original formulation. Among such methods we find \textit{statistical linearization} \cite{Berret20202,Berret20201}, which basically hinges upon the linearization of both the dynamics of the mean and the dynamics of the covariance of the state $x$ around its expected value. More specifically, if only deterministic open-loop controls are leveraged, \eqref{eq:eds} may be replaced by the following system of deterministic differential equations
\begin{equation} \label{eq:syslin}
\begin{cases}
\dot{m}(t) = f(m(t),u(t)), \\
\dot{P}(t) = D_xf(m(t),u(t)) P(t) + P(t) D_xf(m(t),u(t))^{\top} + \G(m(t))\G(m(t))^{\top} ,
\end{cases}
\end{equation}
where $m \in \R^n$ and $P \in M_n(\R)$ denote first-order approximations of the mean and the covariance matrix of the state $x$, respectively. If one is capable of controlling \eqref{eq:syslin} such that additionally $\| P \|$ takes small enough values throughout the control time interval, the quantities $m$ and $P$ will closely resemble the exact mean and covariance matrix of the state $x$, respectively, thus endowing statistical linearization with theoretical guarantees of accurate approximation (indeed, a work in which we carefully analyze the accuracy of this method should appear soon, \textcolor{black}{see \cite{Leparoux2023}}). The main benefit which statistical linearization offers is to reduce minimal variance control, which is a stochastic problem, into an essentially easier problem, which consists of controlling a deterministic system subject to state constraints (on the variable $P$), yielding a numerically cheaper and theoretically guaranteed alternative.

Although this method has already been successfully applied to complex systems, e.g., \cite{Seywald2019,Ridderhof2019}, 
statistical linearization might suffer from controllability issues. Indeed, 
since the variance can be minimized only if the variable $P$ can be controlled directly, at the very least the system \eqref{eq:syslin} must be \emph{accessible}, i.e., the reachable set of \eqref{eq:syslin} must have non-empty interior. It is worth noting that this latter property is fundamental when leveraging statistical linearization for robust optimal control, in that it allows one to locally steer \eqref{eq:syslin} along every directions, a key requirement without which one might not be capable to reduce the value of the cost \cite{Trelat2012,Bonalli2019}. Unfortunately, accessibility might not hold even when considering very simple settings. Indeed, in the case where \eqref{eq:eds} is linear, i.e., the mapping $f$ is linear in $(x,u)$, and the mapping $\G$ is constant, the dynamics of $P$ in \eqref{eq:syslin} is independent from $m$ and $u$, which makes it not accessible: in turn, it is impossible to control the values of $P$ in linear settings.

Although accessibility is difficult to verify and no general exact characterizations exist, 
one may often establish accessibility through sufficient conditions on the rank of a certain Lie algebra (which, in special cases, may additionally yield complete controllability). 
This remark leads us to investigating the following first crucial question for the well-posedness of statistical linearization: are there conditions on the mappings $f$, $\G$ which yield accessibility of \eqref{eq:syslin} by means of Lie algebra-based sufficient conditions? Note that in practice, uncertainty is often difficult to model. 
Therefore, the aforementioned sought conditions must depend on the values of the variable $m$ only (that is, they must not depend on the values of $P$), and in addition, they should be preferably stated in terms of the mapping $f$ uniquely, (that is, they should yield the same properties for different mappings $\G$). 
Once such sufficient conditions are obtained, it is then paramount to understand how restrictive they are. This remark leads us to investigating the following second crucial question for the well-posedness of statistical linearization: do generic systems \eqref{eq:eds} satisfy the aforementioned conditions? The aim of this paper is to answer these two questions.

\textcolor{black}{Before answering the aforementioned questions, let us briefly discuss connections between controllability properties of the original stochastic system and of its statistical linearization, although rigorously establishing such connections is out of the scope of the present work. A priori, the controllability properties of the original stochastic system are not directly related to the controllability properties of its statistical linearization. That said, perhaps under additional conditions one may prove a stochastic control system is ``approximately controllable'' in expectation and in covariance in the case its statistical linearization is controllable. Note that such result would necessarily require to derive estimates of the error between the evolution of the mean and the covariance of the original stochastic system \eqref{eq:eds} and the evolution of the variables $(m,P)$ of its statistical linearization \eqref{eq:syslin}. Importantly, we derive some of such estimates in our upcoming work \cite{Leparoux2023}, partially addressing such matter for some classes of stochastic control systems. In particular, \cite[Corollary 3]{Leparoux2023} and \cite[Proposition 4]{Leparoux2023} yield approximate controllability of the mean of the original stochastic system, although controllability of its covariance still remains an open question.}

The paper is organized as follows. In Section~\ref{sec:Main}, we first summarize existing Lie algebra-based sufficient conditions for accessibility and controllability. Then, we express the Lie brackets generated by the system \eqref{eq:syslin} as functions of the Lie brackets generated by the mapping $f$, thus obtaining sufficient conditions for the accessibility of \eqref{eq:syslin} in terms of the mapping $f$ uniquely. In the particular case of control-affine systems, we additionally obtain conditions for complete controllability. We also study the class of biaffine systems to show that the sufficient conditions we obtain are not necessary in general. Finally, in Section~\ref{sec:transversality} we show that the aforementioned sufficient conditions for accessibility and controllability are generic for stochastic differential equations \eqref{eq:eds} which are control-affine. Specifically, we show that, for some given diffusion $\G$, system \eqref{eq:syslin} is accessible, generically with respect to the drift $f$ and the initial condition $(m_0,P_0)$.

\section{Conditions for accessibility and controllability}
\label{sec:Main}
As we mentioned in the introduction, we consider stochastic control systems of the form~\eqref{eq:eds}, that is
\begin{equation}\nonumber
\mathrm{d}x(t) = f(x(t),u(t)) \; \mathrm{d}t + \G(x(t)) \; \mathrm{d}W_t ,
\end{equation}
where the state $x$ takes values in $\R^n$, the control $u$ takes values in a subset $\U$ of $\R^{\ncontrol}$, and the mappings $f : \R^n \times \R^{\ncontrol} \to \R^n$ and $\G : \R^n \to \R^{n \times \nWiener}$ are assumed to be smooth.
The \emph{statistical linearization} of~\eqref{eq:eds} is the deterministic control system
\begin{equation}\nonumber 
\begin{cases}
\dot{m}(t) = f(m(t),u(t)), \\
\dot{P}(t) = D_xf(m(t),u(t)) P(t) + P(t) D_xf(m(t),u(t))^{\top} + \G(m(t))\G(m(t))^{\top},
\end{cases}
\end{equation}
where $D_xf$ denotes the Jacobian matrix of $f$ with respect to $x$. Here, $m(t) \in \R^n$, $u(t) \in \U$, and $P(t) \in \Sym(n)$, the vector space of $(n\times n)$ symmetric matrices.

Actually, one can show more specific properties of $P(t)$. Indeed, the latter satisfies the differential Lyapunov equation
$$
\dot P = D_xf(m,u)P + P D_x f(m,u)^{\top} + \G(m)\G(m)^{\top} ,
$$
for which it holds
\begin{equation}
\label{eq:intsol}
P(t) = \Phi(t,0)P(0)\Phi(t,0)^{\top} + \int_{0}^t \Phi(t,s) \G(m(s))\G(m(s))^{\top} \Phi(t,s)^{\top} \; \mathrm{d}s , \quad t \ge 0 ,
\end{equation}
where the fundamental matrix $\Phi(t,s) \in \mathbb{R}^{n \times n}$ is defined as the unique solution to
$$
\begin{cases}
\displaystyle \frac{\partial}{\partial t} \Phi(t,s) = D_xf\left(m(t),u(t)\right) \Phi(t,s) , \\[10pt]
\Phi(s,s) = I .
\end{cases}
$$
From \eqref{eq:intsol}, we infer that $P(t)$ is positive definite for every $t\geq 0$ for which $P(t)$ is well-defined, as soon as $P(0)$ is itself positive definite. Motivated by this remark and the fact that the initial value of the covariance of the state $x$ of~\eqref{eq:eds} is selected as initial condition $P(0)$, from now on we assume $P(t)$ is positive definite.

Denote by $\Sym^+(n)$  the set of positive definite $(n\times n)$ symmetric matrices, and consider an open subset $\M$ of $\R^n$ (specifically, in some settings it can be convenient to work with $m$ in $\M$ instead of $m\in \R^n$). We introduce statistical linearization as the control system on $\M \times \Sym^+(n)$ which is given by
\textcolor{black}{\begin{equation}\label{eq:linearization_stat}
  \dot X = F(X,u), \qquad X = \left(\begin{array}{c}
    m \\
    P
  \end{array}\right) \in \M \times \Sym^+(n) , \quad u \in \U,
\end{equation}}
where the smooth mapping $F : \M \times \Sym^+(n) \times \R^{\ncontrol} \to \R^n \times \Sym(n) $ is defined by
\begin{equation}\label{eq:vf_statlin}
F(X, u) = \begin{pmatrix} f(m,u) \\ D_xf(m,u)P + P D_x f(m,u)^{\top} + \G(m)\G(m)^{\top} \end{pmatrix}.
\end{equation}
Our main goal is to study accessibility and controllability properties for this system.

\begin{remark}
\label{re:complete}
As a consequence of~\eqref{eq:intsol}, the trajectories of $F(X,u)$ 
are defined on the same time-intervals in which the trajectories of $f(m,u)$ are. In particular, if for a given $u \in \U$ the vector field $f(\cdot,u)$ is complete on $\R^n$, then the vector field $F(\cdot,u)$ is complete on $\M \times \Sym^+(n)$.
\end{remark}

\subsection{General conditions for accessibility and controllability}
\label{sec:AccConditions}
Let us first recall the main existing accessibility and controllability conditions. We refer to~\cite[Chapter 3]{Jurdjevic1996} for further details. Let
\begin{equation}\label{eq:control_system}
    \dot x = f(x,u), \qquad x \in \M, \  u \in \U \subset \R^{\ncontrol},
\end{equation}
be any (smooth) control system on an open subset $\M$ of $\R^n$. Specifically, this control system is characterized  by the family of (smooth) vector fields
\begin{equation}\nonumber
\f = \lbrace f_u=f(\cdot, u) \ : \  u \in \U \rbrace.
\end{equation}
Two Lie algebras of vector fields are naturally associated with $\f$:
\begin{itemize}
  \item the Lie algebra generated by $\f$,
  \begin{equation}\nonumber
\lie(\f) = \Vect{ \left[ f_{1}, ..., [f_{k-1}, f_{k}] \right]: k\geq 1, f_i \in \f } ;
  \end{equation}
  \item \textcolor{black}{the so called \emph{zero-time ideal}\footnote{\textcolor{black}{The terminology \textit{zero-time ideal} comes from the fact that $\I(\f)$ is associated with the set of chained flows of the form $e^{t_1f_1} \circ \dots \circ e^{t_k f_k}$, $k \in \mathbb{N}$, $f_i \in \f$, where the total time $t_1+\cdots +t_k$ equals zero. On the contrary, $\lie(\f)$ is associated with the set of such elements but for which the total time $t_1+\cdots +t_k$ is let free.}}, defined as}
 \begin{equation}\nonumber
\I(\f) = \Vect{f_{1} - f_{2}, f_3 : f_1, f_2 \in \f, f_3 \in \mathcal{D}(\f) }  \subset \lie(\f),
\end{equation}
where $\mathcal{D}(\f) = \lbrace \left[ f_{1}, ..., [f_{k-1}, f_{k}] \right]: k\geq 2, f_i \in \f \rbrace$. \smallskip
\end{itemize}
From now on, for every $x \in \M$, we denote by $\lie_x(\f)$ (resp.\ $\I_x(\f)$) the set of vectors $f(x) \in \R^n$ such that $f \in \lie(\f)$ (resp.\ $f \in \I(\f)$).

\begin{definition}
We say that $\f$ (or the control system~\eqref{eq:control_system}) satisfies the \emph{accessibility rank condition} at $x \in \M$ if $\dim \lie_x(\f) = n$, and that it satisfies the \emph{accessibility rank condition in fixed time} at $x \in \M$ if $\dim \I_x(\f) = n$.
\end{definition}

\begin{remark}
Note that the accessibility rank condition in fixed time implies the accessibility rank condition since $\I(\f)$ is contained in $\lie(\f)$.
\end{remark}

For every $x \in \M$, let us denote by $\A(x,T)$ the \emph{reachable set} from $x$ in time $T \geq 0$ of~\eqref{eq:control_system}, and by $\A(x)=\bigcup_{T\geq 0} \A(x,T)$ the \emph{reachable set} from $x$. The accessibility rank conditions give sufficient conditions such that reachable sets have nonempty interiors (see also~\cite[Sect.\ 4.5]{Sontag1998} for the necessity of these conditions).

\begin{lemma}[Accessibility properties, Th.\ 2 p.\ 68, and Th.\ 3 p.\ 71 of \cite{Jurdjevic1996}]
\label{lemma:Acc}
\begin{enumerate}
  \item If $\f$ satisfies the accessibility rank condition at every point of $\M$, then $\forall x \in \M$, $\A(x)$ has nonempty interior and the set of interior points is dense in $\A(x)$.
  \item If $\f$ satisfies the accessibility rank condition in fixed time at every point in $\M$, then $\forall x \in \M$, $\A(x,T)$ has nonempty interior for $T>0$ small enough and the set of interior points is dense in $\A(x,T)$. If moreover $\f$ contains at least one complete vector field on $\M$, then the preceding properties of $\A(x,T)$ hold for every $T>0$.
\end{enumerate}
\end{lemma}

The property for a control system of having reachable sets with nonempty interiors is often called \emph{accessibility}, which is a crucial property especially in the context of optimal control. Indeed,  it guarantees that one can generate infinitesimal perturbations in any direction around the endpoint of a trajectory by varying its generating control (either in fixed time if $\mathrm{int} \A(x,T)$ is nonempty, or in free time if $\mathrm{int} \A(x)$ is nonempty), opening up the possibility of performing calculus of variations.\medskip 

Consider now the particular case of control-affine systems, i.e.,
\begin{equation}\label{eq:dynca}
    f(x,u) = f_0(x) + \sum \limits_{i=1}^{\ncontrol} u_i f_i(x), \quad u \in \U=\R^{\ncontrol}.
\end{equation}
The Lie algebras associated with such a system respectively take the forms
\begin{equation}\nonumber
\lie(\f)=\lie \{ f_i  : \ i=0,\dots,\ncontrol \} \quad \hbox{and} \quad \I(\f)=\lie \{ \ad^s f_0 \cdot f_i : \ s\geq 0 , \ i=1,\dots,\ncontrol \} ,
\end{equation}
where we denote $\ad^1 f_0 \cdot f_i = [f_0,f_i]$ and $\ad^s f_0 \cdot f_i = \ad^1 f_0 \cdot (\ad^{s-1} f_0 \cdot f_i)$ for $s \geq 2$. The main benefit offered by these relations is that they allow for expressing the accessibility rank conditions in terms of iterated brackets of the $f_i$'s uniquely. Moreover, a somewhat stronger version of Lemma \ref{lemma:Acc} may be established, enabling the following ease-of-use sufficient conditions for the controllability of control-affine systems.

\begin{lemma}
\label{lemma:Contr}
Assume $\M$ is connected. If for every $x \in \M$ it holds that
\begin{equation}\label{eq:Hormcond}
\dim \lie_x \left( f_1, \dots, f_{\ncontrol} \right) = n,
\end{equation}
then the control-affine system~\eqref{eq:dynca} is \emph{controllable} on $\M$ in free and fixed time, i.e., for every $x \in \M$, $\A(x)=\M$ and $\A(x,T)=\M$ for $T>0$ small enough.
\end{lemma}
\begin{remark}
Note that the Lie algebra $\lie \left( f_1, \dots, f_{\ncontrol} \right)$ is included in $\I(\f)$, therefore \eqref{eq:Hormcond} implies that the accessibility rank condition in fixed time at $x$ holds.
\end{remark}

\subsection{Form of the Lie brackets for statistical linearization}
The accessibility rank conditions require the computation of Lie brackets. We establish here a general expression of the Lie brackets between the vector fields which define the statistical linearization. For this, given a controlled drift $f$ and a diffusion $g$ defining the statistical linearization~\eqref{eq:linearization_stat}, we denote by $\f = \lbrace f(\cdot, u) \ : \  u \in \U \rbrace$ the family of vector fields on $\M$ which is defined by the drift, and by $\F_{\G}= \lbrace F(\cdot, u) \ : \  u \in \U \rbrace$ the family of vector fields~\eqref{eq:vf_statlin} on $\M \times \Sym^+(n)$ which is defined by the statistical linearization. Clearly by definition, for every $F \in \F_{\G}$ there exists $f \in \f$ such that
\textcolor{black}{\begin{equation}\nonumber
F(X) = \vecdef{f(m)}{Df(m)P + PDf(m)^{\top} + \G(m)\G(m)^{\top}}, \quad X = \left(\begin{array}{c}
    m \\
    P
  \end{array}\right) \in \M \times \Sym^+(n) .
\end{equation}}
\begin{lemma} \label{lem:LieBrac}
Consider two vector fields $F_1,F_2$ on $\M \times \Sym^+(n)$ of the form
\begin{equation}\nonumber
F_i(X) = \vecdef{f_i(m)}{Df_i(m)P + PDf_i(m)^{\top} + B_i(m)}, \quad i=1,2,
\end{equation}
where $f_1,f_2$ are vector fields on $\M$ and $B_1,B_2: \M \to  \Sym(n)$ are smooth mappings. Then,
\begin{equation}\nonumber
[F_1, F_2](X) = \begin{pmatrix} [f_1,f_2](m) \\ D[f_1,f_2](m) P + P D[f_1,f_2](m)^{\top} + B_{12}(m)\end{pmatrix},
\end{equation}
where, for every $m \in \M$,
\begin{multline}\nonumber
B_{12}(m) = dB_{2}(m)\cdot f_1(m)  - dB_{1}(m)\cdot f_2(m) + Df_2(m)B_{1}(m) - Df_1(m)B_{2}(m)  \\ + B_{1}(m)Df_2(m)^{\top} - B_{2}(m) Df_1(m)^{\top}
\end{multline} is a symmetric matrix (note that we denote by $Df(m)$ the Jacobian matrix of a vector field $f$, whereas we denote by $dB(m)$ the differential of a map $B$ at $m$).
\end{lemma}

\begin{proof}
By definition, the Lie brackets of $F_1$ and $F_2$ is given by
\[
[F_1,F_2](X) = dF_2(X)\cdot F_1(X) - dF_1(X)\cdot F_2(X) .
\]
Let us denote $Y=(Y_m, Y_P)$ any element $Y \in \R^n \times \Sym(n)$. On the one hand, the first $n$ components of $dF_2(X)\cdot F_1(X)$ satisfy
\begin{multline*}
(dF_2(X)\cdot F_1(X))_m = d(F_2)_m(X)\cdot F_1(X) \\ = \pder{f_2}{m}(m)\cdot (F_{1}(X))_m + \pder{f_{2}}{P}(m)\cdot (F_{1}(X))_P = Df_2(m) f_1(m) ,
\end{multline*}
and therefore we readily obtain that
\begin{equation*}
[F_1,F_2](X)_m = [f_1,f_2](m).
\end{equation*}
On the other hand, the remaining $n^2$ components of $dF_2(X)\cdot F_1(X)$ satisfy
\begin{multline}\nonumber
    (dF_2(X) \cdot F_1(X))_P = \pder{\left(Df_2(m)P + PDf_2(m)^{\top} + B_2(m)\right)}{m}\cdot f_1(m) \\+ \pder{\left(Df_2(m)P + PDf_2(m)^{\top} + B_2(m)\right)}{P}\cdot \left(Df_1(m)P + PDf_1(m)^{\top} + B_1(m)\right).
\end{multline}
Let us develop these two terms separately. For the first term, we may compute
\begin{multline}\nonumber
\pder{\left(Df_2(m)P + PDf_2(m)^{\top}  + B_2(m) \right)}{m}\cdot f_1(m)= \\
D^2f_2(m)\cdot (f_1(m), .)P + P [D^2f_2(m)\cdot (f_1(m), . )]^\top + dB_{2}(m)\cdot f_1(m),
\end{multline}
where, with an abuse of notation, we implicitly identify the linear mapping
$$
D^2f_2(m)\cdot (f_1(m), . ) : h \in \R^n \mapsto D^2f_2(m)\cdot (f_1(m), h) \in \R^n
$$
to an $(n \times n)$ matrix. For the second term, we may compute
\begin{multline}\nonumber
\pder{\left(Df_2(m)P + PDf_2(m)^{\top} + B_2(m) \right)}{P}\cdot \left(Df_1(m)P + PDf_1(m)^{\top} + B_1(m)\right) = \\
 Df_2(m)Df_1(m)P  + Df_2(m)PDf_1(m)^\top + Df_1(m)PDf_2(m)^\top + P  \left[Df_2(m)Df_1(m)\right]^{\top}\\
+  Df_2(m)B_{1}(m) + B_{1}(m)Df_2(m)^{\top}.
\end{multline}
Similar formulas hold for $(dF_1(X)\cdot F_2(X))_P$, which can be readily checked by index exchange. At this step, by leveraging the relation
\begin{multline}\nonumber
D[f_1,f_2](m) = D^2f_2(m)\cdot (f_1(m), .) - D^2f_1(m)\cdot (f_2(m), .) \\
+ Df_2(m)Df_1(m) - Df_1(m)Df_2(m),
\end{multline}
we finally obtain that
\[
[F_1,F_2]_P(X) = D[f_1,f_2](m)P + P D[f_1,f_2](m)^{\top}  + B_{12}(m),
\]
where $B_{12}$ is as stated in the lemma, and the conclusion follows.
\end{proof}

An elementary induction argument on this lemma yields a simple closed form for every element in $\lie (\F_{\G})$ and $\I (\F_{\G})$, respectively.
\begin{corollary}
\label{cor:formofLie}
Every vector field $F$ in $\lie (\F_{\G})$ (resp.\ in $\I (\F_{\G})$) writes as
\begin{equation}\label{eq:Fdans lieF}
F(X) = \vecdef{f(m)}{Df(m)P + PDf(m)^{\top} + B(m)},
\end{equation}
where $f$ belongs to $\lie (\f)$ (resp.\ to $\I (\f)$) and $B: \M \to  \Sym(n)$ is a smooth mapping. Reciprocally, if $f$ belongs to $\lie (\f)$ (resp.\ to $\I (\f)$), then there exists $F$ in $\lie (\F_{\G})$ (resp.\ in $\I (\F_{\G})$) and a smooth mapping $B: \M \to  \Sym(n)$ such that~\eqref{eq:Fdans lieF} holds.

Finally, if $F \in \lie (\F_{0})$ (i.e.\ if $\G=0$), then $B=0$.
\end{corollary}

\begin{remark}
The family $\F_{0}$ plays an important role in the forthcoming algebraic conditions. It is important to note that this family is not directly associated with the statistical linearization of a stochastic differential equation~\eqref{eq:eds}, in that the latter is not defined for $\G=0$. Nevertheless, this family $\F_0$ can be interpreted as the statistical linearization of a deterministic control system subject to imperfect information on the initial condition, which is modelled through a Borel measure (see for instance~\cite{Marigonda2018}).
\end{remark}

\subsection{Algebraic conditions for accessibility}
To simplify the expression of the accessibility rank conditions, we first provide an algebraic result which allows us to express these sought conditions independently from $g$ and $P$. In what follows, we denote $E=\R^n \times \Sym(n)$, a vector space of dimension $N=n + \frac{n(n+1)}{2}$. Recall that $\F_{0}$ denotes the family of the vector fields~\eqref{eq:vf_statlin} associated with $\G=0$. Also, we denote by ${M}_n(\R)$ the space of $(n \times n)$ matrices and by $\Id \in \Sym^+(n)$ the identity matrix.

\begin{lemma}\label{prop:AlgRes}
For every $m \in \M$, the following properties hold:
\begin{enumerate}
  \item $\dim \I_{(m,\Id)} (\F_{0}) = N$ if and only if $\dim \I_{(m,P)} (\F_{0}) = N$ $\forall P \in \Sym^+(n)$;
  \item  $\dim \lie_{(m,\Id)} (\F_{0}) = N$ if and only if $\dim \lie_{(m,P)} (\F_{0}) = N$ $\forall P \in \Sym^+(n)$.
\end{enumerate}
\end{lemma}

\begin{proof}
For any given $P \in \Sym^+(n)$, we introduce the map $\phi_P : \R^n \times {M}_n(\R) \to E$, which is defined by
$$
\phi_P (v,A) = \vecdef{v}{A P + PA^{\top}}.
$$
We claim that
\begin{equation}\label{eq:kerphiP}
\R^n \times M_n(\R) = E \oplus \ker{\phi_P}.
\end{equation}
For this, note first that $\ker{\phi_P} = \left\lbrace 0 \right\rbrace \times \Skew(n) P^{-1}$, where $\Skew(n)$ denotes the set of $(n \times n)$ skew-symmetric matrices. Hence, we are left to prove that
$$
M_n(\R) = \Sym(n) \oplus \Skew(n) P^{-1} ,
$$
or equivalently that
$$
\Sym(n) \cap \Skew(n) P^{-1} = \{0\} .
$$
For this, if we assume $S$ is an element of this intersection, that is $S$ is symmetric and $S = \Lambda P^{-1}$ with $\Lambda \in \Skew(n)$, it must hold that
$$
\Lambda P^{-1} = S = S^{\top} = (\Lambda P^{-1})^{\top} = -P^{-1} \Lambda ,
$$
implying that $P\Lambda + \Lambda P = 0$. This latter Lyapunov equation has no trivial solutions if and only if  $Sp\left\lbrace P\right\rbrace \cap Sp\left\lbrace - P\right\rbrace \neq \emptyset $, which is not the case since $P$ is positive definite. Therefore, $\Lambda = 0$ and thus $S = 0$. We conclude that \eqref{eq:kerphiP} holds true.

At this point, we proceed by first showing the first point of the lemma. For this, fix $m \in \M$ and denote
$$
\widetilde{E} = \left\{ \vecdef{f(m)}{Df(m)} \ : \ f \in \I(\f) \right\}.
$$
Since $\I(\f)$ is a Lie algebra, $\widetilde{E}$ is a linear subspace of $\R^n \times M_n(\R)$. Moreover, for every $P \in \Sym^+(n)$, thanks to Corollary \ref{cor:formofLie} it holds that
$$
\I_{(m,P)} (\F_{0}) = \phi_P(\widetilde{E}) \subset E ,
$$
and hence, thanks to~\eqref{eq:kerphiP} we have that $\dim \I_{(m,P)} (\F_{0}) = N$ if and only if $E \subset \widetilde{E}$. Since this last condition holds true independently from $P$, Point 1.\ of the lemma follows. Point 2.\  can be proved along the same line.
\end{proof}

Let us now study the dimensions of the subspaces $\I_X(\F_{\G})$ and $\lie_X(\F_{\G})$, for $\G \neq 0$. Specifically, the following key density result holds true.

\begin{lemma} \
\label{le:algebraic_condition}
\begin{enumerate}
  \item If $\dim \I_{(m,\Id)} (\F_{0}) = N$ for every $m \in \M$, then there exists an open and dense subset $\Omega$ of $\M \times \Sym^+(n)$, such that $\dim \I_X(\F_{\G}) = N$ for every $X \in\Omega$.
  \item If $\dim \lie_{(m,\Id)} (\F_{0}) = N$ for every $m \in \M$, then there exists an open and dense subset $\widetilde{\Omega}$ of  $\M \times \Sym^+(n)$, such that $\dim \lie_X(\F_{\G}) = N$ for every $X \in \widetilde{\Omega}$.
\end{enumerate}
\end{lemma}

\begin{proof}
We only prove Point 1., the proof of Point 2.\ being identical. Fix $m_0 \in \M$. By assumption, there exist $N$ vector fields $f_1,\dots,f_N$ in $\I(\f)$ such that
\[
\det \nolimits_E \left( \vecdef{f_1(m_0)}{Df_1(m_0) + Df_1(m_0)^{\top}},..., \vecdef{f_N(m_0)}{Df_N(m_0) + Df_N(m_0)^{\top}} \right) \neq 0,
\]
where  $\det \nolimits_E$ denotes the determinant in the vector space $E$.
By continuity, for every $m$ in a neighbourhood $\mathcal{V}_{m_0} \subset \M$ of $m_0$ it holds that
\begin{equation}
\det \nolimits_E \left( \vecdef{f_1(m)}{Df_1(m) + Df_1(m)^{\top}},..., \vecdef{f_N(m)}{Df_N(m) + Df_N(m)^{\top}}  \right) \neq 0.
\label{det}
\end{equation}
At this step, from Corollary~\ref{cor:formofLie}, there exist $F_1, \dots, F_N$ in $\I(\F_{\G})$ such that
$$
F_i(X) = \vecdef{f_i(m)}{Df_i(m)P + PDf_i(m)^{\top} + B_i(m)}, \quad i=1,\dots, N.
$$
Let us introduce the scalar function on $\M \times \Sym^+(n)$, which is defined by
\[
p(X) = \det \nolimits_E \left( F_1(X), ..., F_N (X) \right).
\]
Clearly, $\dim \I_X(\F_{\G}) = N$ as soon as $p(X) \neq 0$. But for any fixed $m \in \mathcal{V}_{m_0} $, the mapping $p(X)$ is a polynomial function of $P$, whose homogeneous term of highest degree is
\[
p^0(X) = \det \nolimits_E \left( \vecdef{f_1(m)}{Df_1(m)P + PDf_1(m)^{\top}},..., \vecdef{f_N(m)}{Df_N(m)P + PDf_N(m)^{\top}}\right).
\]
Therefore, \eqref{det} yields $p^0(m,I) \neq 0$, and thus both $p^0$ and $p$ are not the zero $P$--polynomial. As a consequence, the subset
$$
{\Omega}_{m_0} = \left \lbrace X \in \mathcal{V}_{m_0} \times \Sym^+(n) : \ p(X) \neq 0 \right \rbrace
$$
is open and dense in $\mathcal{V}_{m_0} \times \Sym^+(n)$. By taking the union of the sets ${\Omega}_{m_0}$ for all $m_0 \in \M$, we obtain a dense and open subset $\Omega$ of $\M \times \Sym^+(n)$, on which $\dim \I_X(\F_{\G}) = N$.
\end{proof}

\subsection{Accessibility criterion for statistical linearization}
\label{sec:acc_criterion}
We are now in position to establish sufficient conditions for the accessibility of the statistical linearization. Specifically, the following result comes from a direct application of Lemma~\ref{le:algebraic_condition}.
%

\begin{proposition} \label{prop:acc}
Consider drift $f$ and diffusion $\G$ mappings, and denote by $\f$ and $\F_{\G}$ the families of vector fields associated with $f$ and the corresponding statistical linearization~\eqref{eq:linearization_stat}, respectively. The following conditions hold true:
\begin{enumerate}
    \item if for every $m \in \M$ there holds
 \begin{equation}
 \label{eq:suffcondition1}
     \dim \left\{ \vecdef{f(m)}{Df(m) + Df(m)^{\top}} \ : \ f \in \lie(\f) \right\} = N,
 \end{equation}   
    then $\F_{\G}$ satisfies the accessibility rank condition on an open and dense subset of $\M \times \Sym(n)^+$;
    \item if for every $m \in \M$ there holds
 \begin{equation}
 \label{eq:suffcondition2}
    \dim \left\{ \vecdef{f(m)}{Df(m) + Df(m)^{\top}} \ : \ f \in \I(\f) \right\} = N,
 \end{equation}   
    then $\F_{\G}$ satisfies the accessibility rank condition in fixed time on an open and dense subset of $\M \times \Sym^+(n)$.
\end{enumerate}
\end{proposition}

\begin{remark}
In the case $\G=0$, by Lemma~\ref{prop:AlgRes} the open and dense subsets which result from Proposition \ref{prop:acc} are the whole space $\M \times \Sym^+(n)$, and the requirements of Point 1.\ and 2.\ of the proposition are equivalent to the accessibility rank conditions in free and fixed time, respectively. This equivalence does not hold in general for nonzero $\G$, see Section~\ref{se:bilinear} for a counterexample.
\end{remark}

Importantly, these conditions involve properties of the family of controlled drift mappings $\f$ uniquely, and in particular neither the diffusion $\G$ nor the covariance variable $P$ explicitly affect them. Through these accessibility conditions, one can deduce the structure of the reachable sets in free and fixed time, and thus generic accessibility properties for \eqref{eq:syslin} by means of Lemma~\ref{lemma:Acc}. Note that, thanks to Remark~\ref{re:complete}, the assumption concerning the completeness of the vector fields which is required in Lemma~\ref{lemma:Acc} can be merely tested on the family $\f$.\medskip

Proposition \ref{prop:acc} may be refined in the case of control-affine drift, that is, mappings $f$ of the form
\begin{equation}\nonumber
    f(x,u) = f_0(x) + \sum \limits_{i=1}^{\ncontrol} u_i f_i(x), \quad u \in \U=\R^{\ncontrol} .
\end{equation}
In such case, the statistical linearization takes a control-affine form as well, specifically
\begin{equation}\label{eq:sysca}
\dot{X} = F_{f_0,g}(X) + \sum \limits_{i=1}^{\ncontrol} u_i F_{f_i}(X), \quad u \in \R^{\ncontrol},
\end{equation}
where we introduced the following vector fields on $\M \times \Sym^+(n)$:
\textcolor{black}{$$
\begin{array}{cc}
    F_{f_0,g}(X)=\vecdef{f_0(m)}{Df_0(m)P + P Df_0(m)^{\top} + \G(m)^{\top}\G(m)}, \\[5mm]
    F_{f_i}(X) = \vecdef{f_i(m)}{Df_i(m)P + P Df_i(m)^{\top}}, \qquad X = \left(\begin{array}{c}
    m \\
    P
  \end{array}\right) , \quad  i=1,\dots,{\ncontrol} .
\end{array}
$$}
In addition to the accessibility stated through Proposition~\ref{prop:acc}, from a direct application of Lemma~\ref{lemma:Contr} we obtain stronger controllability conditions for \eqref{eq:sysca}.

\begin{proposition}
\label{prop:accCA}
Assume $\M$ is connected. If for every $m \in \M$ there holds
$$
\dim \left\{ \vecdef{f(m)}{Df(m) + Df(m)^{\top}} \ : \ f \in \lie(f_1,\dots,f_{\ncontrol}) \right\} = N,
$$
then system~\eqref{eq:sysca} is controllable on $\M \times \Sym^+(n)$ in free and fixed time.
\end{proposition}
%
\begin{proof}
It suffices to notice that, by applying Lemma~\ref{prop:AlgRes}-Point 2.\ to the family $\F_0=\{F_{f_1}, \dots, F_{f_{\ncontrol}} \}$, the assumption above implies that \eqref{eq:sysca} satisfies the assumption of Lemma~\ref{lemma:Contr}.
\end{proof}

\subsection{Counterexample: the biaffine case}
\label{se:bilinear}
The aforementioned results give sufficient but not necessary conditions for accessibility. To illustrate that, below we provide an example of a particular class of systems which is generally accessible even though it does not verify the latter sufficient accessibility conditions.

The class we consider is the one associated with constant diffusions $g$ and control-affine drifts $f$ which depend linearly on $x$, i.e.,\ $f_i(x) = A_i x$ with $A_i \in {M}_n(\R)$ for $i= 0,...,m_u$. We say that such an $f$ is a \emph{biaffine system}, the simplest class of nonlinear systems. The statistical linearization takes the control-affine form \eqref{eq:sysca}, where
$$
F_{f_0,g}(X)=\vecdef{A_0m }{A_0P + P A_0^{\top} + \G\G^{\top}}, \quad
    F_{f_i}(X) = \vecdef{A_im }{A_iP + P A_i^{\top}}, \ i=1,\dots,\ncontrol .
$$
For simplicity we set $F_0=F_{f_0,g}$ and $F_i=F_{f_i}$. Lemma \ref{lem:LieBrac} yields the following forms of the Lie brackets,
\begin{equation}
\label{eq:brackets_biaffine}
[F_i , F_j](X) = \begin{pmatrix} \left[ A_i , A_j \right] m \\ \left[ A_i , A_j \right] P + P \left[ A_i , A_j \right]^{\top} + B_{ij} \end{pmatrix} , \quad i , j = 0 , \dots , \ncontrol , 
\end{equation}
where the symmetric matrices $B_{ij} \in \Sym(n)$ are given by
\begin{equation}
\label{eq:Bij}
B_{0j}=A_j g g^\top + g g^\top A_j^\top \qquad \hbox{and}  \qquad B_{ij} = 0 \qquad \hbox{for } i , j = 1, \dots , \ncontrol.
\end{equation}

We first show that, in such specific case, the accessibility conditions given by Proposition \ref{prop:acc} do no longer hold true.

\begin{proposition}\label{prop: biliNAcc}
Assume that the drift $f$ is a biaffine system. Then, at every $m \in \R^n$, the sufficient conditions~\eqref{eq:suffcondition1} and~\eqref{eq:suffcondition2} in Proposition \ref{prop:acc} are not satisfied. In particular, if in addition $\G = 0$, then the statistical linearization  satisfy none of the accessibility rank conditions.
\end{proposition}

\begin{proof}

 As a consequence of \eqref{eq:brackets_biaffine}, condition~\eqref{eq:suffcondition1} at  $m \in \R^n$ writes as
\begin{equation*} 
    \dim \left\{\vecdef{A m}{A + A^{\top}} : \ A \in \lie ( A_0 , A_1 , \dots, A_{\ncontrol} ) \right\} = N .
\end{equation*}
Fix $m \in \R^n$, $P \in \Sym^+(n)$, and introduce the linear map $\psi_{(m,P)} : {M}_n(\R) \to \R^n \times \Sym(n)$ defined by
$$
\psi_{(m,P)} (A) = \vecdef{A m}{AP +PA^{\top}}.
$$
If condition~\eqref{eq:suffcondition1} holds at $m$, then necessarily $\mathrm{rank} \psi_{(m,I)} =N$. Let us prove that this latter condition is not actually true.

For this, obviously we have that
$$
\ker \psi_{(m,I)} = \left \lbrace A \in \Skew(n) : \ A m = 0 \right \rbrace .
$$
If $m=0$, then $\ker \psi_{(m,I)}=\Skew(n)$ is of dimension $n(n-1)/2$, and thus, by the rank-nullity theorem, $\mathrm{rank} \psi_{(m,I)} = n^2 - n(n-1)/2= n(n+1)/2$ is smaller than $N$.
If $m \neq 0$, then we can find an orthogonal basis $\{ e_1 , \dots , e_n \}$ of $\R^n$ where $e_1 = m$. In that basis, the element $A_{i,j}$ in the $i{th}$ row and the $j{th}$ column of $A$ satisfies
$$
A_{i,j} =  e_i^{\top} A e_j = -e_j^{\top} A e_i = - A_{j,i} ,
$$
yielding $A_{1,j} = A_{i,1} = 0$. As a consequence, in these coordinates, we have
$$
\ker \psi_{(m,I)} = \left\lbrace A = \begin{pmatrix} 0&0\\0 & \bar{A} \end{pmatrix} \in {M}_n(\R) : \ \bar{A} \in \Skew(n-1) \right \rbrace .
$$
Thus $\ker \psi_{(m,I)}$ has the same dimension as $\Skew(n-1)$. Hence $\mathrm{rank} \psi_{(m,I)} = n^2 - (n-2)(n-1)/2$ is equal to $N-1$ and we deduce that the condition~\eqref{eq:suffcondition1} is never satisfied. This proof can be easily replicated to show condition~\eqref{eq:suffcondition2}.
\end{proof}

Proposition \ref{prop: biliNAcc} shows that our sufficient conditions for accessibility can not be leveraged to prove that biaffine systems are generally accessible. Nevertheless, it is still possible to formulate other sufficient general conditions for the accessibility of biaffine systems. In particular, the following result holds true.

\begin{proposition}
Assume that the biaffine system satisfies the followings:
\begin{enumerate}
    \item[(i)] $\dim \lie (A_1 , \dots, A_{\ncontrol})  = n^2$;
    \item[(ii)] there exists $i \in \{1,\dots,\ncontrol\}$ such that $B_{0i}=A_i g g^\top + g g^\top A_i^\top$ is nonzero.
\end{enumerate}
Then  $\F_\G$ satisfies the accessibility rank condition in fixed time on an open and dense subset of $ \R^n \times \Sym^+(n)$. 
\end{proposition}

\begin{remark}
Conditions \emph{(i)} and \emph{(ii)} are satisfied by generic matrices $A_1,\dots, A_{\ncontrol}$ when $\ncontrol >1$ and $g \neq 0$, respectively. Hence the statistical linearizations of a ``very large'' class of biaffine systems satisfy the accessibility rank conditions, though they do not satisfy the sufficient conditions which are provided in Proposition \ref{prop:acc}.
\end{remark}

\begin{proof}
Note first that, as a consequence of \eqref{eq:brackets_biaffine}, every vector field in $\I(\F_g)$ writes as
$$
F(X) = \vecdef{A m}{AP + PA^{\top} + B}, \qquad A \in M_n(\R), \quad B \in \Sym(n),
$$
and thus $F(X)$ depends affinely on $X=(m,P)$. As a result,  the accessibility rank condition in fixed time $\dim \I_X(\F_g)=N$ holds either nowhere, or for points $X$ which belong to an open and dense subset of $ \R^n \times \Sym^+(n)$ (to see this, use the argument we provide at the end of the proof of Lemma \ref{le:algebraic_condition}). Therefore, it is enough to prove that $\dim \I_X(\F_g)=N$ holds for any specific point $X$.

For this, by Hypothesis \emph{(ii)} the symmetric matrix $B_{0i}$ is nonzero, thus we can choose $\bar{m} \in \R^n$ and $\bar{P} \in \Sym^+(n)$ such that $\bar{m}^\top\bar{P}^{-1} B_{0i} \bar{P}^{-1}\bar{m} \neq 0$. Set $\bar{X}=(\bar{m},\bar{P})$.

Now, note that, on the one hand the set $\I_{\bar{X}}(\F_g)$ contains the vector
$$
[F_0,F_i](\bar{X}) = \vecdef{A_{0i} \bar{m}}{A_{0i}\bar{P} + \bar{P}A_{0i}^{\top} + B_{0i}}, \quad \hbox{where } A_{0i} = \left[ A_0 , A_i \right] ,
$$
whereas, on the other hand, it contains the linear subset
$$
W = \left \lbrace F(\bar{X}) \ : \ F \in \lie(F_1,\dots,F_{\ncontrol}) \right \rbrace .
$$
Thanks to \eqref{eq:brackets_biaffine} and \eqref{eq:Bij}, every $F(\bar{X}) \in W$ writes as
$$
F(\bar{X}) = \vecdef{A \bar{m}}{A\bar{P} + \bar{P}A^{\top}}, \qquad \hbox{with } A \in \lie( A_1 , \dots, A_{\ncontrol} ).
$$
It then results from Hypothesis \emph{(i)} that $W$ is the image of the map $\psi_{\bar{X}}$ which is defined in the proof of Proposition~\ref{prop: biliNAcc}. Moreover, a straightforward adaptation of the latter proof implies that $\mathrm{rank}\psi_{\bar{X}}= N-1$, and therefore that $\dim W = N-1$.

Let us prove that $[F_0,F_i](\bar{X})$ does not belong to $W$. Indeed, if an element $(v,Q) \in \R^n \times \Sym(n)$ belongs to $W$, then there exists a matrix $A \in M_n(\R)$ such that $v=A \bar{m}$ and $Q= A \bar{P}+ \bar{P}A^{\top}$. Equivalently, there exists  $\Omega \in \Skew(n)$ such that  $v = (\frac{1}{2} Q + \Omega) \bar{P}^{-1}\bar{m}$. By multiplying the latter expression by $\bar{m}^{\top}\bar{P}^{-1}$, we obtain that
$$
\alpha(v,Q) :=\bar{m}^{\top}\bar{P}^{-1}(v - \frac{1}{2} Q \bar{P}^{-1}\bar{m}) = 0.
$$
Let us compute this expression in the case $(v,Q)=[F_0,F_i](\bar{X})$. We have that
\begin{multline*}
\alpha([F_0,F_i](\bar{X})) = \bar{m}^{\top}\bar{P}^{-1} \left(A_{0i} \bar{m} - \frac{1}{2} (A_{0i}\bar{P} + \bar{P}A_{0i}^{\top} + B_{0i}) \bar{P}^{-1}\bar{m} \right) \\
= \frac{1}{2} \bar{m}^{\top} \bar{P}^{-1} \left( (A_{0i}\bar{P} - \bar{P}A_{0i}^{\top}) - B_{0i} \right)\bar{P}^{-1} \bar{m} = -\frac{1}{2} \bar{m}^\top \bar{P}^{-1} B_{0i} \bar{P}^{-1} \bar{m}.
\end{multline*}
But the choice of $\bar{X}$ guarantees $\alpha([F_0,F_i](\bar{X})) \neq 0$, which implies $[F_0,F_i](\bar{X}) \not\in W$. 

In summary, we have shown that  $W \oplus \R [F_0,F_i](\bar{X})$ is a vector space of dimension $N$ included in $\I_{\bar{X}}(\F_g)$.  We deduce that $\dim \I_{\bar{X}}(\F_g)=N$, concluding the proof.
\end{proof}

\section{Generic accessibility and controllability}
\label{sec:transversality}
In this section, we conclude our study by establishing generic accessibility and controllability properties for \eqref{eq:syslin}. Specifically, we show that the various sufficient rank conditions  which we established in Section \ref{sec:acc_criterion} hold \textit{generically} (in the sense of Thom's transversality) with respect to the drift $f$, and eventually with respect to the initial conditions $m(0)$ and $P(0)$. 

From now on, we consider control-affine drift mappings uniquely, i.e.,
$$
f(x,u) = f_0(x) + \sum^{\ncontrol}_{i=1} u_i f_i(x) , \quad u \in \U=\R^{\ncontrol} ,
$$
and we thus adopt the notation we introduced through \eqref{eq:sysca}. For an open subset $\M$ of $\R^n$ and $\ell , r \in \mathbb{N}$, we introduce the following quantities (see, e.g., \cite{Hirsch1994}):
\begin{ItemizE}
    \item $\pi^{\ell}_r : J^{\ell}(\M,\mathbb{R}^r) \rightarrow \M$ is the vector bundle of $\ell$-jets of maps in $C^{\infty}(\M,\mathbb{R}^r)$;
    \item for every $h \in C^{\infty}(\M,\mathbb{R}^r)$, $j^{\ell}_r(h) : \M \rightarrow J^{\ell}(\M,\mathbb{R}^r)$ is the jet map associated with $\pi^{\ell}_r$; in particular, when clear from the context, for the sake of conciseness we denote $j^{\ell}(h)_y = j^{\ell}(h)(y) \in J^{\ell}(\M,\mathbb{R}^r)$, for every $h \in C^{\infty}(\M,\mathbb{R}^r)$ and every $y \in \M$.
\end{ItemizE}
We recall that the family of sets
$$
\Big\{ h \in C^{\infty}(\M,\mathbb{R}^r) : \ j^{\ell}(h)(\M) \subset \mathcal{V} \Big\} , \quad \textnormal{$\mathcal{V}$ open subset of $J^{\ell}(\M,\mathbb{R}^r)$} ,
$$
generates the (weak) Whitney topology of $C^{\infty}(\M,\mathbb{R}^r)$. Our main results are as follows.


\textcolor{black}{\begin{theorem}[Accessibility with Respect to Generic Drift] \label{theo:Acc}\label{theo:Acc+}
There exists an open and dense (w.r.t.\ the Whitney topology) subset $\mathcal{O}_{\textnormal{Acc}} \subset C^{\infty}(\M,\mathbb{R}^{n \times ({\ncontrol} + 1)})$, such that for every $(f_0,f_1,\dots,f_{\ncontrol}) \in \mathcal{O}_{\textnormal{Acc}}$, the family of vector fields
$$
\f = \left\{ f_0 + \sum_{i=1}^{\ncontrol}u_i f_i \ : \ u \in \mathcal{U} \right\}
$$
satisfies condition~\eqref{eq:suffcondition1} for every $m \in \M$. In addition, there exists an open and dense (w.r.t.\ the Whitney topology) subset $\mathcal{O}_{\textnormal{Acc}_+} \subset C^{\infty}(\M,\mathbb{R}^{n \times ({\ncontrol} + 1)})$, such that for every $(f_0,f_1,\dots,f_{\ncontrol}) \in \mathcal{O}_{\textnormal{Acc}_+}$, the aforedefined family of vector fields $\f$ satisfies condition~\eqref{eq:suffcondition2} for every $m \in \M$.
\end{theorem}}



Together with Proposition \ref{prop:acc}, Theorem \ref{theo:Acc+} readily yields the following corollary.

\begin{corollary}[Accessibility with Respect to Generic Drift and Initial Datum] \label{coroll:AccTotal}
There exists an open and dense (w.r.t.\ the Whitney topology) subset $\mathcal{O}_{\textnormal{Acc}_+} \subset C^{\infty}(\M,\mathbb{R}^{n \times ({\ncontrol} + 1)})$, such that for every $(f_0,f_1,\dots,f_{\ncontrol}) \in \mathcal{O}_{\textnormal{Acc}_+}$ and every diffusion $g \in C^{\infty}(\M,\mathbb{R}^{n \times \nWiener})$, 
the family of vector fields
$$
\F = \left\{ F_{f_0,g} + \sum_{i=1}^{\ncontrol}u_i F_{f_i} \ : \ u \in \mathcal{U} \right\}
$$
satisfies the accessibility rank condition in fixed time on an open and dense subset $\Omega_g \subset \M \times \Sym^+(n)$ which depends on $g$ uniquely.
\end{corollary}

\begin{theorem}[Controllability with Respect to Generic Drift] \label{theo:Contr}\ 

Let ${\ncontrol} \ge 2$. There exists an open and dense (w.r.t.\  the Whitney topology) subset $\mathcal{O}_{\textnormal{Contr}} \subset C^{\infty}(\M,\mathbb{R}^{n \times {\ncontrol}})$, such that for every $(f_1,\dots,f_{\ncontrol}) \in \mathcal{O}_{\textnormal{Contr}}$, every vector field $f_0 \in C^{\infty}(\M,\mathbb{R}^n)$, and every diffusion $g \in C^{\infty}(\M,\mathbb{R}^{n \times \nWiener})$, 
the family of vector fields
$$
\F = \left\{ F_{f_0,g} + \sum_{i=1}^{\ncontrol}u_i F_{f_i} \ : \ u \in \mathcal{U} \right\}
$$
is controllable on $\M \times \Sym^+(n)$ in free and fixed time.
\end{theorem}

Through the rank conditions provided above, one can deduce the structure of the reachable sets in free and fixed time of the statistical linearization, and thus generic accessibility and controllability properties for \eqref{eq:syslin}, by means of Lemma~\ref{lemma:Acc}--\ref{lemma:Contr}.

The rest of the section is devoted to proving these results. In particular, we prove \textcolor{black}{the first part of} Theorem \ref{theo:Acc} for ${\ncontrol} = 1$ (the cases ${\ncontrol} \ge 2$ trivially result from the latter), given that the proofs of \textcolor{black}{Theorem \ref{theo:Contr} and of the second part of Theorem \ref{theo:Acc+}} follow the same line, though they are less pedagogical. Moreover, since all our proofs develop along local analysis, without loss of generality we implicitly assume $\M = \R^n$ in what follows. 

\subsection{Fundamental concepts from transversality theory}

Before jumping into the proof of Theorem \ref{theo:Acc}, we recall some useful concepts from algebraic geometry and transversality theory which are key to achieve the aforementioned proof, for which we mainly refer to \cite{Whitney1957,Bochnak1998,Martinet1970,Goresky1988}. In particular, for the sake of clarity in the exposition, we gathered the classical concepts of algebraic geometry we make use of in our proofs in Appendix \ref{sec:appendixWhitney}, so that in this subsection we may focus on transversality theory.

Since we set up ${\ncontrol} = 1$ and $\M = \R^n$, throughout the following sections, for any given $\ell \in \mathbb{N}$ we consider the vector bundle of $\ell$-jets $J^{\ell}(\mathbb{R}^n,\mathbb{R}^{2 n})$, for which we implicitly adopt the following identification
\begin{equation} \label{eq:Identification}
    J^{\ell}(\mathbb{R}^n,\mathbb{R}^{2 n}) \cong \mathbb{R}^{3 n + 2 n \sum^{\ell}_{i=1} {\scriptsize \left(\begin{array}{c}
    n + i - 1 \\
    i
\end{array}\right)}} ,
\end{equation}
which is achieved through some (global) trivialization mapping. In addition, if $(U,\varphi)$ is any local chart of $\mathbb{R}^n$ at $m \in \mathbb{R}^n$, for every $f_0 , f_1 \in C^{\infty}(\mathbb{R}^n,\mathbb{R}^n)$ the coordinates $j^{\ell}_{\varphi}(f_0,f_1)_m$ of the jet $j^{\ell}(f_0,f_1)_m \in J^{\ell}(\mathbb{R}^n,\mathbb{R}^{2 n})$ with respect to $(U,\varphi)$ are denoted by
\begin{align*}
j^{\ell}_{\varphi}( f_0 , f_1 )_m = \bigg( \varphi(m) , \ &\varphi^{f_0}_0 , \varphi^{f_1}_0 , \\
&(\varphi^{f_0}_{1;\{i_1\}})_{0 \le i_1 \le n} , (\varphi^{f_1}_{1;\{i_1\}})_{0 \le i_1 \le n} , \\
&\dots\dots\dots\dots\dots\dots\dots\dots\dots\dots \\
&(\varphi^{f_0}_{\ell;\{i_1,\dots,i_{\ell}\}})_{0 \le i_1 \le \dots \le i_{\ell} \le n} ,  (\varphi^{f_1}_{\ell;\{i_1,\dots,i_{\ell}\}})_{0 \le i_1 \le \dots \le i_{\ell} \le n} \bigg) ,
\end{align*}
where $(\varphi^{f_0}_j)_{0 \le i_1 \le \dots \le i_j \le n} , (\varphi^{f_1}_j)_{0 \le i_1 \le \dots \le i_j \le n} \in \mathbb{R}^{n {\scriptsize \left(\begin{array}{c}
    n + j - 1 \\
    j
\end{array}\right)}}$, for every $1 \le j \le \ell$.

\vspace{5pt}


Through the identification \eqref{eq:Identification}, 
we may introduce semi-algebraic sets $\mathcal{SA} \subset J^{\ell}(\mathbb{R}^n,\mathbb{R}^{2 n})$ 
(see Appendix \ref{sec:appendixWhitney}), and thus we can make use of the following key transversality result to prove Theorem \ref{theo:Acc}.

\begin{theorem}[Transversality Theorem, Proposition p. 38
of \cite{Goresky1988}] \label{theo:Transversality}
Let $\mathcal{SA} \subset J^{\ell}(\mathbb{R}^n,\mathbb{R}^{2 n})$ be a closed semi-algebraic set which satisfies
$$
\codim \mathcal{SA} \ge n + 1 .
$$
Then, the subset
$$
\mathcal{O}_{\mathcal{SA}} = \Big\{ (f_0,f_1) \in C^{\infty}(\mathbb{R}^n,\mathbb{R}^{2 n}) : \ j^{\ell}(f_0,f_1) \cap \mathcal{SA} = \emptyset \Big\}
$$
is open and dense (w.r.t. the Whitney topology) in $C^{\infty}(\mathbb{R}^n,\mathbb{R}^{2 n})$.
\end{theorem}

Theorem \ref{theo:Transversality} results from combining Theorem \ref{theo:Whitney} with Thom's stratified transversality theorem (for details, see \cite[``Th\'eor\`eme 6.1'' and ``Remarque'' at page 175]{Martinet1970}).

\subsection{Proof of Theorem \ref{theo:Acc}}

Our proof of Theorem \ref{theo:Acc} is inspired by the proofs in \cite{Lobry1972,Bonnard1997,Chitour2006}, which consist of showing that specific \textit{bad sets} have arbitrarily large codimension. We introduce appropriate, non-trivial modifications to these proofs to address submersion singularities which are due to possibly non-transverse components. In particular, we introduce four bad sets $\mathcal{B}^1_{\ell}$, $\mathcal{B}^2_{\ell}$, $\mathcal{B}^3_{\ell}$, and $\mathcal{B}^4_{\ell}$ (see below). Although we could prove our result without explicitly considering $\mathcal{B}^1_{\ell}$, we introduced this latter on purpose to show how the classical transversality cunning in \cite{Lobry1972,Bonnard1997,Chitour2006}, which can be leveraged to handle $\mathcal{B}^1_{\ell}$, must be revisited to manipulate $\mathcal{B}^2_{\ell}$, $\mathcal{B}^3_{\ell}$, and $\mathcal{B}^4_{\ell}$.

\subsubsection{Definition of the bad sets} Through a slight abuse of notation which is possible thanks to the identification \eqref{eq:Identification} 
, for $r \in \mathbb{N} \cap [0,\ell-2]$ the mappings
$$
j^{\ell}(f_0,f_1)_m \in J^{\ell}(\mathbb{R}^n,\mathbb{R}^{2 n}) \mapsto \ad^r F_{f_0,0} \cdot F_{f_1}(m,\Id) \in \mathbb{R}^N ,
$$
$$
j^{\ell}(f_0,f_1)_m \in J^{\ell}(\mathbb{R}^n,\mathbb{R}^{2 n}) \mapsto \ad^r \big[ F_{f_0,0} , F_{f_1} \big] \cdot F_{f_1}(m,\Id) \in \mathbb{R}^N ,
$$
where $\ad^r$ is the $r$-iterated Lie bracket, are well-defined vector-valued polynomials. This remark, together with the identification \eqref{eq:Identification} 
, allows us to define the following bad sets, which are semi-algebraic sets of $J^{\ell}(\mathbb{R}^n,\mathbb{R}^{2 n})$ 
(for any $\ell \in \mathbb{N}$ large enough):
\begin{align*}
    \mathcal{B}^1_{\ell} = \bigg\{ &j^{\ell}(f_0,f_1)_m \in J^{\ell}(\mathbb{R}^n,\mathbb{R}^{2 n}) : \ f_0(m) \neq 0 , \\
    &\rank \Big[ \ad^1 F_{f_0,0} \cdot F_{f_1}(m,\Id) \ \Big| \ \ad^3 F_{f_0,0} \cdot F_{f_1}(m,\Id) \ \Big| \ \dots \\
    &\hspace{29ex} \dots \ \Big| \ \ad^{2 \lfloor \ell / 2 \rfloor - 1} F_{f_0,0} \cdot F_{f_1}(m,\Id) \Big] < N \bigg\} ,
\end{align*}
\begin{align*}
    \mathcal{B}^2_{\ell} = \bigg\{ &j^{\ell}(f_0,f_1)_m \in J^{\ell}(\mathbb{R}^n,\mathbb{R}^{2 n}) : \ f_0(m) = 0 , \ \ad^1 f_0 \cdot f_1(m) \neq 0 , \\
    &\rank \Big[ \ad^1 \big[ F_{f_0,0} , F_{f_1} \big] \cdot F_{f_1}(m,\Id) \ \Big| \ \ad^3 \big[ F_{f_0,0} , F_{f_1} \big] \cdot F_{f_1}(m,\Id) \ \Big| \ \dots \\
    &\hspace{18ex} \dots \ \Big| \ \ad^{2 \lfloor (\ell - 1) / 2 \rfloor - 1} \big[ F_{f_0,0} , F_{f_1} \big] \cdot F_{f_1}(m,\Id) \Big] < N \bigg\} ,
\end{align*}
\begin{align*}
    \mathcal{B}^3_{\ell} = \bigg\{ &j^{\ell}(f_0,f_1)_m \in J^{\ell}(\mathbb{R}^n,\mathbb{R}^{2 n}) : \ f_0(m) = \ad^1 f_0 \cdot f_1(m) = 0 , \quad f_1(m) \neq 0 \bigg\} ,
\end{align*}
\begin{align*}
    \mathcal{B}^4_{\ell} = \bigg\{ &j^{\ell}(f_0,f_1)_m \in J^{\ell}(\mathbb{R}^n,\mathbb{R}^{2 n}) : \ f_0(m) = f_1(m) = 0 \bigg\} ,
\end{align*}

Thanks to \eqref{eq:DimClosure}
, we finally have that the bad set
$$
\mathcal{B}_{\ell} = \overline{ \mathcal{B}^1_{\ell} \cup \mathcal{B}^2_{\ell} \cup \mathcal{B}^3_{\ell} \cup \mathcal{B}^4_{\ell}} \subset J^{\ell}(\mathbb{R}^n,\mathbb{R}^{2 n})
$$
is a semi-algebraic set whose codimension satisfies
\begin{equation} \label{eq:TestCodim}
    \codim \mathcal{B}_{\ell} \ge \underset{i=1,2,3,4}{\min} \ \codim \mathcal{B}^i_{\ell} .
\end{equation}

\subsubsection{Existence of an open and dense subset}

Assume that
\begin{equation} \label{eq:ReqCodim}
    \underset{i=1,2,3,4}{\min} \ \codim \mathcal{B}^i_{\ell} \ge n + 1 , \quad \textnormal{for some} \ \ell \in \mathbb{N} .
\end{equation}
Therefore, \eqref{eq:TestCodim} yields $\codim \mathcal{B}_{\ell} \ge n + 1$ and we can apply Theorem \ref{theo:Transversality} to infer the existence of an open and dense (w.r.t.\ the Whitney topology) subset
$$
\mathcal{O}_{\textnormal{Acc}} = \Big\{ (f_0,f_1) \in C^{\infty}(\mathbb{R}^n,\mathbb{R}^{2 n}) : \ j^{\ell}(f_0,f_1) \cap \mathcal{B}_{\ell} = \emptyset \Big\} \subset C^{\infty}(\mathbb{R}^n,\mathbb{R}^{2 n}) .
$$
In particular, straightforward computations by contradiction show that
$$
(f_0,f_1) \in \mathcal{O}_{\textnormal{Acc}} \ \Longrightarrow \ \dim \lie_m \left( \vecdef{f}{Df + Df^{\top}} \ : \ f \in \{ f_0 , f_1 \} \right) = N .
$$
Therefore, the proof of Theorem \ref{theo:Acc} is achieved once \eqref{eq:ReqCodim} is proved.

\begin{remark}
Theorem \ref{theo:Acc+} may be proved similarly by replacing $\mathcal{B}^1_{\ell}$ with
\begin{align*}
    \bigg\{ &j^{\ell}(f_0,f_1)_m \in J^{\ell}(\mathbb{R}^n,\mathbb{R}^{2 n}) : \ f_0(m) \neq 0 , \ \ad^1 f_0 \cdot f_1(m) \neq 0 , \\
    &\rank \Big[ \ad^1 \big[ F_{f_0,0} , F_{f_1} \big] \cdot F_{f_1}(m,\Id) \ \Big| \ \ad^3 \big[ F_{f_0,0} , F_{f_1} \big] \cdot F_{f_1}(m,\Id) \ \Big| \ \dots \\
    &\hspace{18ex} \dots \ \Big| \ \ad^{2 \lfloor (\ell - 1) / 2 \rfloor - 1} \big[ F_{f_0,0} , F_{f_1} \big] \cdot F_{f_1}(m,\Id) \Big] < N \bigg\} ,
\end{align*}
The required modifications are straightforward, and thus we avoid to report them.
\end{remark}

\subsubsection{Codimension computations}

The rest of the manuscript is devoted to prove \eqref{eq:ReqCodim}. We provide computations for $\mathcal{B}^1_{\ell}$ and $\mathcal{B}^2_{\ell}$ uniquely, given that the computations for $\mathcal{B}^3_{\ell}$ and $\mathcal{B}^4_{\ell}$ essentially follow the same line, though they are easier.

We fix $\ell \in \mathbb{N}$ large enough (to be selected later). For $r \in \mathbb{N}$, $m \in \mathbb{R}^n$, and $h_0 , \dots , h_r \in C^{\infty}(\mathbb{R}^n,\mathbb{R}^n)$, we denote by
$$
\big( D^r h_0 \big)_{i_1,\dots,i_r}(m) = \frac{\partial^r h_0}{\partial x_{i_1} \dots \partial x_{i_r}}(m) \in \mathbb{R}^n
$$
the tensor-valued mapping given by the $r$th derivative of $h_0$ with respect to the multi-index $(i_1,\dots,i_r)$, as well as the matrix-valued mapping
$$
\big( (D^{r+1} h_0) h_1 \dots h_r \big)_{i j}(m) = \sum^n_{i_1,\dots,i_r=1} \big( D^{r+1} h_0 \big)^i_{j,i_1,\dots,i_r}(m) h^{i_1}_1(m) \dots h^{i_r}_r(m) \in \mathbb{R}^{n \times n} .
$$
Also, we denote the symmetric tensorization of any element $\theta_r \in \mathbb{R}^{n {\scriptsize \begin{pmatrix}
    n + r - 1 \\
    r
\end{pmatrix}}}$ by
$$
S(\theta_r) : \underbrace{\mathbb{R}^n \times \dots \times \mathbb{R}^n}_{r-\textnormal{times}} \rightarrow \mathbb{R}^n ,
$$
whereas we denote the symmetrization of a matrix $A \in \mathbb{R}^n$ by $\Sym(A) = A + A^{\top}$. For every $r \in \mathbb{N}$ and $\rho \in [0,r]$, we denote by $L^{\rho}_r$ the Stiefel manifold of linear mappings
$$
L : \mathbb{R}^r \rightarrow \mathbb{R}^N
$$
such that $\rank L = \rho$, for which it holds that
\begin{equation} \label{eq:Stiefel}
    \codim L^{\rho}_r = \left( N - \rho \right) (r - \rho) .
\end{equation}
Finally, for every $r \in \mathbb{N}$, we introduce the notation $\bm{1}_r = (1,\dots,1) \in \mathbb{R}^r$.

\vspace{10pt}

\textbf{Computations for $\bm{\mathcal{B}^1_{\ell}}$.}

\vspace{5pt}

We introduce appropriate coordinates, i.e., the mapping $\chi$ defined below, under which it is easy to check that appropriate mappings are submersions, yielding the desired result. Note that the mapping $\chi$ has been first introduced in \cite[Section 5]{Bonnard1997}. Unfortunately, due to some Lie bracket singularities, the technique developed in \cite{Bonnard1997} does not work for $\mathcal{B}^2_{\ell}$, and we will need to introduce new singularity-free coordinates.

\vspace{5pt}

Straightforward computations by induction yield the following key result.

\begin{proposition} \label{Prop:FirstBrackets}
For 
$m \in \mathbb{R}^n$ and $f_0 \in C^{\infty}(\mathbb{R}^n,\mathbb{R}^n)$ such that $f_0(m) \neq 0$, the following equivalence holds true for $f_1 \in C^{\infty}(\mathbb{R}^n,\mathbb{R}^n)$ and $r \in \mathbb{N}$:
\begin{align*}
    \ad^r F_{f_0,0} \cdot F_{f_1}(m,\Id) = \ &\left(\begin{array}{c}
        A_{r}\big( D^{i+1} f_0 , D^i f_1 , \ 0 \le i \le r-1 \big) \\[10pt]
        B_{r}\big( D^{i+1} f_0 , D^i f_1 , \ 0 \le i \le r \big)
    \end{array}\right) \\[10pt]
    &+ \left(\begin{array}{c}
        (D^r f_1) \underbrace{f_0 \dots f_0}_{r-\textnormal{times}} \\
        \Sym\Big( (D^{r+1} f_1) \underbrace{f_0 \dots f_0}_{r-\textnormal{times}} \Big)
    \end{array}\right) ,
\end{align*}
where $A_r$ and $B_r$ are polynomial expressions of 
the derivatives of $f_0$ and $f_1$.
\end{proposition}

Consider the canonical global chart $(\mathbb{R}^n,Id)$ of $\mathbb{R}^n$, given by the identity map $Id$. Since every 
$j^{\ell}(f_0,f_1)_m \in \mathcal{B}^1_{\ell}$ satisfies $Id^{f_0}_0 \neq 0$
, we may select constant vectors $v_2 , \dots , v_n \in \mathbb{R}^n$ such that the tuple of vectors
$$
\{ v_1(p) , v_2(p) , \dots , v_n(p) \} = \{ Id^{f_0}_0 , v_2 , \dots , v_n \} , \quad \textnormal{for} \  p = j^{\ell}(f_0,f_1)_m \in \mathcal{V}^1_{\ell} ,
$$
spans $\mathbb{R}^n$ in some open subset 
$\mathcal{V}^1_{\ell} \subset J^{\ell}(\mathbb{R}^n,\mathbb{R}^{2 n})$. 
In particular, the mapping
\begin{align*}
    \chi\big( &j^{\ell}(f_0,f_1)_m \big) = \bigg( \chi^m_0 = m , \\
    &\chi^{f_0}_0 = Id^{f_0}_0 , \quad \chi^{f_1}_0 = Id^{f_1}_0 , \\
    &\Big( \chi^{f_0}_{1;\{i_1\}} \Big)_{1 \le i_1 \le n} = \Big( S(Id^{f_0}_1) v_{i_1} \Big)_{1 \le i_1 \le n} , \\
    &\Big( \chi^{f_1}_{1;\{i_1\}} \Big)_{1 \le i_1 \le n} = \Big( S(Id^{f_1}_1) v_{i_1} \Big)_{1 \le i_1 \le n} , \\
    &\dots\dots\dots\dots\dots\dots\dots\dots\dots\dots\dots\dots\dots \\
    &\Big( \chi^{f_0}_{\ell;\{i_1,\dots,i_{\ell}\}} \Big)_{1 \le i_1 \le \dots \le i_{\ell} \le n} = \Big( S(Id^{f_0}_{\ell}) v_{i_1} \dots v_{i_{\ell}} \Big)_{1 \le i_1 \le \dots \le i_{\ell} \le n} , \\
    &\Big( \chi^{f_1}_{\ell;\{i_1,\dots,i_{\ell}\}} \Big)_{1 \le i_1 \le \dots \le i_{\ell} \le n} = \Big( S(Id^{f_1}_{\ell}) v_{i_1} \dots v_{i_{\ell}} \Big)_{1 \le i_1 \le \dots \le i_{\ell} \le n} \bigg) ,
\end{align*}
which is smoothly well-defined in $\mathcal{V}^1_{\ell}$, is a local chart for 
$J^{\ell}(\mathbb{R}^n,\mathbb{R}^{2 n})$. For 
$j^{\ell}(f_0,f_1)_m \in \mathcal{V}^1_{\ell}$, let $M$ denote the matrix which corresponds to the smooth change of coordinates
$$
e^j = \sum^n_{i=1} M_{ij}(\theta^{f_0}_0) v_i , \quad M(\theta^{f_0}_0) = \big( M_{ij}(\theta^{f_0}_0) \big)_{1 \le i , j \le n} \in GL(n) ,
$$
where $e^j$ denotes the $j$th column of $\Id \in \Sym^+(n)$. 
Therefore, for $r \in \mathbb{N}$ we compute
\begin{align*}
    (D^{r+1} \ &f_1) \underbrace{f_0 \dots f_0}_{r-\textnormal{times}} = \left( S(Id^{f_1}_{r+1}) \underbrace{Id^{f_0}_0 \dots Id^{f_0}_0}_{r-\textnormal{times}} e^1 \Big| \ \dots \ \Big| S(Id^{f_1}_{r+1}) \underbrace{Id^{f_0}_0 \dots Id^{f_0}_0}_{r-\textnormal{times}} e^n \right) \\[10pt]
    &= \left( S(Id^{f_1}_{r+1}) \underbrace{v_1 \dots v_1}_{r-\textnormal{times}} v_1 \Big| S(Id^{f_1}_{r+1}) \underbrace{v_1 \dots v_1}_{r-\textnormal{times}} v_2 \Big| \ \dots \ \Big| S(Id^{f_0}_{r+1}) \underbrace{v_1 \dots v_1}_{r-\textnormal{times}} v_n \right) M \\[10pt]
    &= \Big( \chi^{f_1}_{r+1;\bm{1}_{r+1}} \Big| \chi^{f_1}_{r+1;\{\bm{1}_r,2\}} \Big| \ \dots \ \Big| \chi^{f_1}_{r+1;\{\bm{1}_r,n\}} \Big) M ,
\end{align*}
which together with the equivalence of Proposition \ref{Prop:FirstBrackets} finally yields
\begin{align} \label{eq:TransformedBracketFirst}
    \ad^r F_{f_0,0} \cdot \ &F_{f_1}(m,\Id) = \left(\begin{array}{c}
        A_{r}\big( \chi^{f_0}_{i+1} , \chi^{f_1}_i , \ 0 \le i \le r-1 \big) \\[10pt]
        B_{r}\big( \chi^{f_0}_{i+1} , \chi^{f_1}_i , \ 0 \le i \le r \big)
    \end{array}\right) \\[10pt]
    &+ \left(\begin{array}{c}
        \chi^{f_1}_r \\
        \Sym\bigg( \Big( \chi^{f_1}_{r+1;\bm{1}_{r+1}} \Big| \chi^{f_1}_{r+1;\{\bm{1}_r,2\}} \Big| \ \dots \ \Big| \chi^{f_1}_{r+1;\{\bm{1}_r,n\}} \Big) M(\chi^{f_0}_0) \bigg)
    \end{array}\right) , \nonumber
\end{align}
for every 
$j^{\ell}(f_0,f_1)_m \in \mathcal{V}^1_{\ell}$ and every $r \in \mathbb{N} \cap [1,\ell-1]$.

At this step, up to appropriate identifications, the triangular form of \eqref{eq:TransformedBracketFirst} allows us to infer that, for every $\ell \in \mathbb{N}$ such that $2 \lfloor \ell / 2 \rfloor - 1 \ge N$, the mappings

\vspace{5pt}

\begin{align*}
    \mathcal{L}^1_{\ell} &= \Big[ \ad^1 F_{f_0,0} \cdot F_{f_1}(m,\Id) \ \Big| \ \ad^3 F_{f_0,0} \cdot F_{f_1}(m,\Id) \ \Big| \ \dots \\
    &\hspace{20ex} \dots \ \Big| \ \ad^{2 \lfloor \ell / 2 \rfloor - 1} F_{f_0,0} \cdot F_{f_1}(m,\Id) \Big] : \ \mathcal{V}^1_{\ell} \rightarrow \mathbb{R}^{N \times (2 \lfloor \ell / 2 \rfloor - 1)}
\end{align*}

\vspace{5pt}

\noindent are submersions. In particular, we infer that $(\mathcal{L}^1_{\ell})^{-1}(L^{\rho}_{2 \lfloor \ell / 2 \rfloor - 1}) \subset \mathcal{V}^1_{\ell}$ are submanifolds for every $\rho \in [0,N]$, which due to \eqref{eq:Stiefel} satisfy
\begin{align} \label{eq:Stiefel1}
    \codim (\mathcal{L}^1_{\ell})^{-1}(L^{\rho}_{2 \lfloor \ell / 2 \rfloor - 1}) &= \codim L^{\rho}_{2 \lfloor \ell / 2 \rfloor - 1} \\
    &= \left( N - \rho \right) \Big( 2 \lfloor \ell / 2 \rfloor - 1 - \rho \Big) . \nonumber
\end{align}
Since
$$
\mathcal{B}^1_{\ell} \cap \mathcal{V}^1_{\ell} \subset \bigcup^{N - 1}_{\rho = 1} (\mathcal{L}^1_{\ell})^{-1}(L^{\rho}_{2 \lfloor \ell / 2 \rfloor - 1}) ,
$$
thanks to \eqref{eq:Stiefel1} we finally obtain that
\begin{align*}
    \codim \mathcal{B}^1_{\ell} \cap \mathcal{V}^1_{\ell} &\ge \codim (\mathcal{L}^1_{\ell})^{-1}\left( L^{N - 1}_{2 \lfloor \ell / 2 \rfloor - 1} \right) \\
    &= 2 \lfloor \ell / 2 \rfloor - n - \frac{n(n+1)}{2} \ge n + 1 ,
\end{align*}
as soon as we select $\ell \in \mathbb{N}$ large enough. We conclude that \eqref{eq:ReqCodim} holds true for $\mathcal{B}^1_{\ell}$.

\vspace{10pt}

\textbf{Computations for $\bm{\mathcal{B}^2_{\ell}}$.}

\vspace{5pt}

Unfortunately, the technique developed in the previous section does not work for $\mathcal{B}^2_{\ell}$, in that the vector field $f_0$ is identically zero. Thus, we introduce new singularity-free coordinates, i.e., the mapping $\psi$ defined below, which allow us to easily show surjectivity of some differentials, yielding the desired result also in this singular case.

\vspace{5pt}

We start proving the existence of new singularity-free coordinates enabling the following key equivalence between high-order Lie brackets.

\begin{proposition} \label{Prop:SecondBrackets}
For every 
$m \in \mathbb{R}^n$ and every $f_0 , f_1 \in C^{\infty}(\mathbb{R}^n,\mathbb{R}^n)$ satisfying $\ad^1 f_0 \cdot f_1(m) \neq 0$, there exists a local chart $(U,\varphi)$ of $\mathbb{R}^n$ at $m$ which smoothly depends on the coefficients of both $f_0$ and $f_1$, and such that for every $r \in \mathbb{N}$, the following equivalence holds true:
\begin{align*}
    \ad^r \big[ F_{f_0,0} , F_{f_1} \big] \cdot \ &F_{f_1}(m,\Id) = \left(\begin{array}{c}
        0 \\[10pt]
        \displaystyle \Sym\left( (D^2 \varphi^{-1}) \frac{\partial^r f_1}{\partial^r x^{\varphi}_1} D\varphi \right)
    \end{array}\right) \\[10pt]
    &+ \left(\begin{array}{c}
        \displaystyle D \varphi^{-1} \frac{\partial^r f_1}{\partial^r x^{\varphi}_1} \\[10pt]
        \displaystyle \Sym\left( D \varphi^{-1} \left( \frac{\partial^{r+1} f_1}{\partial^{r+1} x^{\varphi}_1} \ \bigg| \ \frac{\partial^{r+1} f_1}{\partial^r x^{\varphi}_1 x^{\varphi}_2} \ \bigg| \ \dots \ \bigg| \ \frac{\partial^{r+1} f_1}{\partial^r x^{\varphi}_1 x^{\varphi}_n} \right) D \varphi \right)
    \end{array}\right) .
\end{align*}
\end{proposition}

\begin{proof}[Proof of Proposition \ref{Prop:SecondBrackets}]
We first start by recalling a classic representation result, in a general abstract framework. Let $M$ be a $n$-dimensional manifold, and $V \in C^{\infty}(M;TM)$ be a smooth vector field, i.e., a section of the tangent bundle $\pi_M : TM \to M$. Fix $p \in M$, and let $(U_1,\varphi)$ and $(U_2,\psi)$ be local charts of $M$ at $p$, so that
$$
\tilde \varphi\left( p , \sum^n_{i=1} v_i \partial^{\varphi}_i(p) \right) = (\varphi(p) , v) , \quad \tilde \psi\left( p , \sum^n_{i=1} w_i \partial^{\psi}_i(p) \right) = (\varphi(p) , w)
$$
are local charts for $TM$ with domains $\pi^{-1}_M(U_1)$ and $\pi^{-1}_M(U_2)$, respectively. If
$$
\sum^n_{i=1} V^{\varphi}_i(p) \partial^{\varphi}_i(p) = V(p) = \sum^n_{i=1} V^{\psi}_i(p) \partial^{\psi}_i(p) ,
$$
the following rules for the change of coordinates of $V$ and its differential, hold:
\begin{equation} \label{eq_change1}
    V^{\psi}(p) = D\big( \psi \circ \varphi^{-1} \big)(p) V^{\varphi}(p) ,
\end{equation}
\begin{align} \label{eq_change2}
    D\big( \tilde \psi \circ V \circ \psi^{-1} \big)(p) &= \\
    &= D\big( \tilde \psi \circ \tilde \varphi^{-1} \big)(p,V(p)) D\big( \tilde \varphi \circ V \circ \varphi^{-1} \big)(p) D\big( \varphi \circ \psi^{-1} \big)(p) , \nonumber
\end{align}
with
$$
D\big( \tilde \psi \circ \tilde \varphi^{-1} \big)(p,V(p)) = \left(\begin{array}{cc}
    D\big( \psi \circ \varphi^{-1} \big)(p) & 0 \\[10pt]
    D^2\big( \psi \circ \varphi^{-1} \big) V^{\varphi}(p) & D\big( \psi \circ \varphi^{-1} \big)(p)
\end{array}\right) \in GL(2 n) .
$$

At this step, since $\ad^1 f_0 \cdot f_1(m) \neq 0$ we may find a local chart $(U,\varphi)$ of $\mathbb{R}^n$ at $m$, which smoothly depends on the coefficients of both $f_0$ and $f_1$ (indeed, this chart is a smooth mapping of the flow of the vector field $\ad^1 f_0 \cdot f_1$, which is smooth with respect to the coefficients of both $f_0$ and $f_1$), and such that
$$
[ f_0 , f_1 ] = \partial^{\varphi}_1 , \quad \textnormal{locally in} \ U .
$$
In these coordinates, through straightforward computations, we obtain
$$
( \ad^r [ f_0 , f_1 ] \cdot f_1 )^{\varphi} = \frac{\partial^r f_1}{\partial^r x^{\varphi}_1} , \quad \textnormal{locally in} \ U ,
$$
and therefore, by selecting $M = \mathbb{R}^n$, $V = \ad^r [ f_0 , f_1 ] \cdot f_1$, $(U_1,\varphi) = (U,\varphi)$, and $(U_2,\psi) = (\mathbb{R}^n,Id)$, the relations \eqref{eq_change1} and \eqref{eq_change2} respectively yield
$$
\ad^r [ f_0 , f_1 ] \cdot f_1 = D \varphi^{-1} \frac{\partial^r f_1}{\partial^r x^{\varphi}_1} , \quad \textnormal{locally in} \ U ,
$$
\begin{align*}
    &D\big( \ad^r [ f_0 , f_1 ] \cdot f_1 \big) = (D^2 \varphi^{-1}) \frac{\partial^r f_1}{\partial^r x^{\varphi}_1} D\varphi \\
    &+ D \varphi^{-1} \left( \frac{\partial^{r+1} f_1}{\partial^{r+1} x^{\varphi}_1} \ \bigg| \ \frac{\partial^{r+1} f_1}{\partial^r x^{\varphi}_1 x^{\varphi}_2} \ \bigg| \ \dots \ \bigg| \ \frac{\partial^{r+1} f_1}{\partial^r x^{\varphi}_1 x^{\varphi}_n} \right) D \varphi , \quad \textnormal{locally in} \ U ,
\end{align*}
and the conclusion readily follows from Lemma \ref{lem:LieBrac}.
\end{proof}

Let us introduce the subset
$$
\mathcal{S}^2_{\ell} = \bigg\{ j^{\ell}(f_0,f_1)_m \in J^{\ell}(\mathbb{R}^n,\mathbb{R}^{2 n}) : \ f_0(m) = 0 , \ \ad^1 f_0 \cdot f_1(m) \neq 0 \bigg\} .
$$
Since the mapping
$$
\mathcal{F}^2_{\ell} : \ J^{\ell}(\mathbb{R}^n,\mathbb{R}^{2 n}) \rightarrow \mathbb{R}^n : \ j^{\ell}(f_0,f_1)_m \mapsto f_0(m)
$$
is a submersion, we infer that
$$
\mathcal{S}^2_{\ell} = ( \mathcal{F}^2_{\ell} )^{-1}(0) \setminus \Big\{ j^{\ell}(f_0,f_1)_m \in J^{\ell}(\mathbb{R}^n,\mathbb{R}^{2 n}) : \ \ad^1 f_0 \cdot f_1(m) = 0 \Big\}
$$
is an embedded submanifold of 
$J^{\ell}(\mathbb{R}^n,\mathbb{R}^{2 n})$ with 
$\dim \mathcal{S}^2_{\ell} = \dim J^{\ell}(\mathbb{R}^n,\mathbb{R}^{2 n}) - n$, which contains $\mathcal{B}^2_{\ell}$. By fixing $m \in \mathbb{R}^n$ and letting $(U,\varphi)$ be the local chart of $\mathbb{R}^n$ at $m$ given by Proposition \ref{Prop:SecondBrackets}, which smoothly depends on the coefficients of both $f_0$ and $f_1$, we define the open set 
$\mathcal{V}^2_{\ell} = \Pi^{-1}_{\ell}(U) \subset J^{\ell}(\mathbb{R}^n,\mathbb{R}^{2 n})$. Therefore, the mapping
\begin{align*}
    \psi\big( j^{\ell}(f_0,f_1)_m \big) = &\bigg( \psi^m_0 = \varphi(m) , \\
    &\psi^{f_0}_0 = D \varphi^{-1}(m) \varphi^{f_0}_0 , \quad \psi^{f_1}_0 = D \varphi^{-1}(m) \varphi^{f_1}_0 , \\
    &\Big( \psi^{f_0}_{1;\{i_1\}} \Big)_{1 \le i_1 \le n} = \Big( D \varphi^{-1}(m) \varphi^{f_0}_{1;\{i_1\}} \Big)_{1 \le i_1 \le n} , \\
    &\Big( \psi^{f_1}_{1;\{i_1\}} \Big)_{1 \le i_1 \le n} = \Big( D \varphi^{-1}(m) \varphi^{f_1}_{1;\{i_1\}} \Big)_{1 \le i_1 \le n} , \\
    &\dots\dots\dots\dots\dots\dots\dots\dots\dots\dots\dots\dots\dots\dots\dots \\
    &\Big( \psi^{f_0}_{\ell;\{i_1,\dots,i_{\ell}\}} \Big)_{1 \le i_1 \le \dots \le i_{\ell} \le n} = \Big( D \varphi^{-1}(m) \varphi^{f_0}_{\ell;\{i_1,\dots,i_{\ell}\}} \Big)_{1 \le i_1 \le \dots \le i_{\ell} \le n} , \\
    &\Big( \psi^{f_1}_{\ell;\{i_1,\dots,i_{\ell}\}} \Big)_{1 \le i_1 \le \dots \le i_{\ell} \le n} = \Big( D \varphi^{-1}(m) \varphi^{f_1}_{\ell;\{i_1,\dots,i_{\ell}\}} \Big)_{1 \le i_1 \le \dots \le i_{\ell} \le n} \bigg) ,
\end{align*}
that is smoothly well-defined in $\mathcal{V}^2_{\ell}$, is a local chart for 
$J^{\ell}(\mathbb{R}^n,\mathbb{R}^{2 n})$ which is adapted to $\mathcal{S}^2_{\ell}$. Finally, evaluating the equivalence of Proposition \ref{Prop:SecondBrackets} through this chart yields
\begin{align} \label{eq:TransformedBracketSecond}
    \ad^r \big[ &F_{f_0,0} , F_{f_1} \big] \cdot F_{f_1}(m,\Id) = \nonumber \\
    &= \left(\begin{array}{c}
        0 \\[10pt]
        \displaystyle \Sym\left( (D^2 \varphi^{-1})(\psi^m_0) D \varphi\big( \varphi^{-1}(\psi^m_0) \big) \psi^{f_1}_{r;\bm{1}_r} D \varphi\big( \varphi^{-1}(\psi^m_0) \big) \right)
    \end{array}\right) \\[10pt]
    &+ \left(\begin{array}{c}
        \displaystyle \psi^{f_1}_{r;\bm{1}_r} \\[10pt]
        \displaystyle \Sym\bigg( \Big( \psi^{f_1}_{r+1;\bm{1}_{r+1}} \ \Big| \ \psi^{f_1}_{r+1;\{\bm{1}_r,2\}} \ \Big| \ \dots \ \Big| \ \psi^{f_1}_{r+1;\{\bm{1}_r,n\}} \Big) D \varphi\big( \varphi^{-1}(\psi^m_0) \big) \bigg)
    \end{array}\right) , \nonumber
\end{align}
for every 
$j^{\ell}(f_0,f_1)_m \in \mathcal{S}^2_{\ell} \cap \mathcal{V}^2_{\ell}$ and every $r \in \mathbb{N} \cap [0,\ell-2]$.

At this step, up to appropriate identifications, the triangular form of \eqref{eq:TransformedBracketSecond} allows us to infer that, for every $\ell \in \mathbb{N}$ such that $2 \lfloor (\ell - 1) / 2 \rfloor - 1 \ge N$, the mappings
\begin{align*}
    \mathcal{L}^2_{\ell} = \Big[ \ad^1 & \big[ F_{f_0,0} , F_{f_1} \big] \cdot F_{f_1}(X) \ \Big| \ \ad^3 \big[ F_{f_0,0} , F_{f_1} \big] \cdot F_{f_1}(X) \ \Big| \ \dots \\
    &\dots \ \Big| \ \ad^{2 \lfloor (\ell - 1) / 2 \rfloor - 1} \big[ F_{f_0,0} , F_{f_1} \big] \cdot F_{f_1}(X) \Big] : \ \mathcal{S}^2_{\ell} \cap \mathcal{V}^2_{\ell} \rightarrow \mathbb{R}^{N \times (2 \lfloor (\ell - 1) / 2 \rfloor - 1)}
\end{align*}
are submersions. In particular, we infer that $(\mathcal{L}^2_{\ell})^{-1}(L^{\rho}_{2 \lfloor (\ell - 1) / 2 \rfloor - 1}) \subset \mathcal{S}^2_{\ell} \cap \mathcal{V}^2_{\ell}$ are submanifolds for every $\rho \in [0,N]$, which due to \eqref{eq:Stiefel} satisfy
\begin{align} \label{eq:Stiefel2}
    \codim (\mathcal{L}^2_{\ell})^{-1}(L^{\rho}_{2 \lfloor (\ell - 1) / 2 \rfloor - 1}) &= \codim L^{\rho}_{2 \lfloor (\ell - 1) / 2 \rfloor - 1} \\
    &= \left( N - \rho \right) \Big( 2 \lfloor (\ell - 1) / 2 \rfloor - 1 - \rho \Big) . \nonumber
\end{align}
Since
$$
\mathcal{B}^2_{\ell} \cap \mathcal{V}^2_{\ell} \subset \bigcup^{N - 1}_{\rho = 1} (\mathcal{L}^2_{\ell})^{-1}(L^{\rho}_{2 \lfloor (\ell - 1) / 2 \rfloor - 1}) ,
$$
thanks to \eqref{eq:Stiefel2} we finally obtain that
\begin{align*}
    \codim \mathcal{B}^2_{\ell} \cap \mathcal{V}^2_{\ell} &\ge \codim (\mathcal{L}^2_{\ell})^{-1}\left( L^{N - 1}_{2 \lfloor (\ell - 1) / 2 \rfloor - 1} \right) \\
    &= 2 \lfloor (\ell - 1) / 2 \rfloor - n - \frac{n(n+1)}{2} \ge n + 1 ,
\end{align*}
as soon as we select $\ell \in \mathbb{N}$ large enough. We conclude that \eqref{eq:ReqCodim} holds true for $\mathcal{B}^2_{\ell}$.

This latter argument concludes the proof of Theorem \ref{theo:Acc}.

\appendix\section{Useful concepts from algebraic geometry}
\label{sec:appendixWhitney}
Let $\ell , r \in \mathbb{N}$. A semi-algebraic set $\mathcal{SA} \subseteq \mathbb{R}^r$ is a subset of the form
\begin{equation} \label{eq:SemiAlg}
    \mathcal{SA} = \bigcup^{\ell}_{i = 1} \Big\{ y \in \mathbb{R}^r : \ p^i_1(y) = \dots = p^i_{a_i}(y) = 0 , \ q^i_1(y) < 0 , \dots , q^i_{b_i}(y) < 0 \Big\} ,
\end{equation}
for polynomials $p^i_j , q^i_j : \mathbb{R}^r \rightarrow \mathbb{R}$. Semi-algebraic sets enjoy the following properties:
\begin{ItemizE}
    \item finite unions of semi-algebraic sets are semi-algebraic sets;
    \item closures (w.r.t. the natural topology) of semi-algebraic sets are semi-algebraic sets.
\end{ItemizE}
A simple induction argument on \eqref{eq:SemiAlg} shows that semi-algebraic sets are disjoint unions of open subsets of algebraic sets. Thus, Whitney's stratification theorem readily yields the following representation result for semi-algebraic sets.

\begin{theorem}[Whitney's stratification for semi-algebraic sets {\cite[Th.~2, p.~546]{Whitney1957}}] \label{theo:Whitney}
Let $r \in \mathbb{N}$ and $\mathcal{SA} \subseteq \mathbb{R}^r$ be a semi-algebraic set. There exist $\ell \in \mathbb{N}$ and embedded submanifolds $\mathcal{M}_1 , \dots , \mathcal{M}_{\ell} \subseteq \mathbb{R}^r$, which are mutually disjoint and such that
$$
\mathcal{SA} = \bigcup^{\ell}_{i = 1} \mathcal{M}_i \quad \textnormal{and} \quad \dim \mathcal{M}_{i+1} < \dim \mathcal{M}_i , \ i = 1 , \dots , \ell-1 .
$$
\end{theorem}

Thanks to Theorem \ref{theo:Whitney}, given $r \in \mathbb{N}$ and a semi-algebraic set $\mathcal{SA} \subseteq \mathbb{R}^r$, we may define the dimension and the codimension of $\mathcal{SA}$ as
$$
\dim \mathcal{SA} = \underset{i = 1,\dots,\ell}{\max} \dim \mathcal{M}_i \quad \textnormal{and} \quad \codim \mathcal{SA} = r - \dim \mathcal{SA} .
$$
Importantly, the topological closure does not modify the dimension, nor the codimension of a semi-algebraic set. Indeed, as a matter of example let $r \in \mathbb{N}$ and
$$
\mathcal{SA} = \big\{ y \in \mathbb{R}^r : \ p(y) = 0 , \ q(y) < 0 \big\} ,
$$
for polynomials $p , q : \mathbb{R}^r \rightarrow \mathbb{R}$. By definition, we must have
\begin{align*}
    \dim \overline{\mathcal{SA}} &\le \max \Big\{ \dim \mathcal{SA} \ , \ \dim \big\{ y \in \mathbb{R}^r : \ p(y) = q(y) = 0 \big\} \Big\} \\
    &\le \max \Big\{ \dim \mathcal{SA} \ , \ \dim \big\{ y \in \mathbb{R}^r : \ p(y) = 0 \big\} \Big\} = \dim \mathcal{SA} ,
\end{align*}
where we infer the last equality from the fact that $\mathcal{SA}$ is open in $p^{-1}(0)$, and thus
\begin{equation} \label{eq:DimClosure}
    \dim \overline{\mathcal{SA}} = \dim \mathcal{SA} \quad \textnormal{and} \quad \codim \overline{\mathcal{SA}} = \codim \mathcal{SA} .
\end{equation}
The extension of this property to general semi-algebraic sets follows by iterating the previous argument through classical induction.


\bibliographystyle{siamplain}
\bibliography{main.bib}
\end{document}